\documentclass[12pt]{article}
\usepackage{latexsym}
\usepackage{amssymb}
\parskip=\medskipamount
\newtheorem{definition}{Definition}[section]
\newtheorem{lemma}[definition]{Lemma}
\newtheorem{theorem}[definition]{Theorem}
\newtheorem{proposition}[definition]{Proposition}

\newtheorem{remark}[definition]{Remark}

\newtheorem{examples}[definition]{Examples}

 1   
\font\ddpp=msbm10  scaled \magstep 1  
\newfont{\msi}{msbm8 scaled \magstephalf}  

\def\gpd{\rightrightarrows}
\def\gr{G\gpd M}
\def\grr{G\times \R\gpd M\times \R}
\def\s{\alpha}
\def\t{\beta}
\def\mult{\sigma :G\to \R}
\def\Jacobi{\Lambda , E}
\def\sostJ{\# _{(\Lambda ,E)}}
\def\sostP{\#_\Lambda}
\def\alg{(A,\lcf \, , \, \rcf ,\rho )}
\def\estalg{(\lcf \, , \, \rcf ,\rho )}
\def\QED{\hskip0.1em\hfill\null\ \null\nobreak\hfill
\kern3pt\lower1.8pt\vbox{\hrule\hbox
{\vrule\kern1pt\vbox{\kern1.7pt \hbox{$\scriptstyle
QED$}\kern0.2pt}\kern1pt\vrule}\hrule}}

\let\izq=\overleftarrow
\let\der=\overrightarrow
\def\ssli{{\stackrel{\leftarrow}{X}}}
\def\sslib{{\stackrel{\leftarrow}{Y}}}
\def\expsn{e^{-\sigma (g)}}
\def\exps{e^{\sigma (g)}}
\def\tLambda{\tilde{\Lambda}}
\def\tssli{{\stackrel{\leftarrow}{\tilde{X}}}}

\def\inv{^{-1}}

\def\fract{\frac{\partial}{\partial t}}

\def\Rp{\hbox{\msi R}}           
\def\R{\hbox{\ddpp R}}               

\def\lcf{\lbrack\! \lbrack}
\def\rcf{\rbrack\! \rbrack}

\newcommand\prueba {\mbox{{\em Proof: }}}

\begin{document}
\baselineskip=.55cm
\title{{\bf JACOBI GROUPOIDS  AND GENERALIZED LIE BIALGEBROIDS}}
\author{David Iglesias-Ponte, Juan C. Marrero
\\ {\small\it Departamento de Matem\'atica Fundamental, Facultad de
Matem\'aticas,}\\[-8pt] {\small\it Universidad de la Laguna, La
Laguna,} \\[-8pt] {\small\it Tenerife, Canary Islands,
Spain,}\\[-8pt] {\small\it E-mail: diglesia@ull.es,
jcmarrer@ull.es} }
\date{\empty}

\maketitle \baselineskip=.45cm
\begin{abstract}
{\small \baselineskip .45cm Jacobi groupoids are introduced as a
generalization of Poisson and contact groupoids and it is proved
that generalized Lie bialgebroids are the infinitesimal invariants
of Jacobi groupoids. Several examples are discussed.}
\end{abstract}
\begin{quote}

\end{quote}

\baselineskip=.5cm
\section{Introduction}
A Poisson groupoid is a Lie groupoid $\gr$ with a Poisson
structure $\Lambda$ for which the graph of the partial
multiplication is a coisotropic submanifold in the Poisson
manifold $(G\times G\times G,\Lambda \oplus \Lambda \oplus
-\Lambda )$ (see \cite{We2}). If $(\gr ,\Lambda )$ is a Poisson
groupoid then there exists a Poisson structure on $M$ such that
the source projection $\s :G\to M$ is a Poisson morphism.
Moreover, if $AG$ is the Lie algebroid of $G$ then the dual bundle
$A^\ast G$ to $AG$ itself has a Lie algebroid structure. Poisson
groupoids were introduced by Weinstein \cite{We2} as a
generalization of both Poisson Lie groups and the symplectic
groupoids which arise in the integration of arbitrary Poisson
manifolds. A canonical example of symplectic groupoid is the
cotangent bundle $T^\ast G$ of an arbitrary Lie groupoid $\gr$. In
this case, the base space is $A^\ast G$ and the Poisson structure
on $A^\ast G$ is just the linear Poisson structure induced by the
Lie algebroid $AG$ (see \cite{CDW}).

In \cite{MX}, Mackenzie and Xu proved that a Lie groupoid $\gr$
endowed with a Poisson structure $\Lambda$ is a Poisson groupoid
if and only if the bundle map $\sostP:T^\ast G\to TG$ is a
morphism between the cotangent groupoid $T^\ast G\gpd A^\ast G$
and the tangent groupoid $TG\gpd TM$. This characterization was
used in order to prove that Lie bialgebroids are the infinitesimal
invariants of Poisson groupoids, that is, if $(\gr ,\Lambda )$ is
a Poisson groupoid then $(AG,A^\ast G)$ is a Lie bialgebroid and,
conversely, a Lie bialgebroid structure on the Lie algebroid of a
(suitably simply connected) Lie groupoid can be integrated to a
Poisson groupoid structure \cite{LX,MX,MX2} (these results can be
applied to obtain a new proof of a theorem of Karasaev \cite{Ka}
and Weinstein \cite{We1} about the relation between symplectic
groupoids and their base Poisson manifolds). We remark that in
\cite{CL}, Crainic and Fernandes have given the precise
obstructions to integrate an arbitrary Lie algebroid to a Lie
groupoid.

On the other hand, a contact groupoid $(\gr ,\eta ,\sigma )$ is
Lie groupoid $\gr$ endowed with a contact 1-form $\eta \in \Omega
^1 (G)$ and a multiplicative function $\sigma \in C^\infty (G,\R)$
such that $$\eta _{(gh)}(X_g \oplus _{TG} Y_h )=\eta _g(X_g)+\exps
\eta _h(Y_h),\mbox{ for }(X_g,Y_h)\in TG^{(2)},$$ where $\oplus
_{TG}$ is the partial multiplication in the tangent Lie groupoid
$TG\gpd TM$ (see \cite{D1,D2,KS}). Contact groupoids can be
considered as the odd-dimensional counterpart of symplectic
groupoids and they have applications in the prequantization of
Poisson manifolds and in the integration of local Lie algebras
associated to rank one vector bundles (see \cite{D1,D2}). In this
case, the base space $M$ carries an induced Jacobi structure such
that the pair $(\s ,e^\sigma )$ is a conformal Jacobi morphism.
Moreover, the presence of the multiplicative function $\sigma$
induces a 1-cocycle $\phi _0\in \Gamma (A^\ast G)$ in the Lie
algebroid cohomology of $AG$. We note that the relation between
Jacobi structures and Lie algebroids with 1-cocycles has been
recently explored in \cite{IM0} by the authors. More precisely, we
have obtained that a Lie algebroid structure on a vector bundle
$A\to M$ and a 1-cocycle $\phi _0\in \Gamma (A^\ast )$, a
generalized Lie algebroid in our terminology, induce a Jacobi
structure $(\Lambda _{(A^\ast ,\phi _0)},E_{(A^\ast ,\phi _0)})$
on $A^\ast$ satisfying some linearity conditions. In addition,
using the differential calculus on Lie algebroids in the presence
of a 1-cocycle, it has been introduced in \cite{IM} (see also
\cite{GM}) the notion of a generalized Lie bialgebroid in such a
way that a Jacobi manifold has associated a canonical generalized
Lie bialgebroid. A generalized Lie bialgebroid is a pair $((A,\phi
_0),(A^\ast ,X_0))$, where $(A,\phi _0)$ and $(A^\ast ,X_0)$ are
generalized Lie algebroids, such that the Lie algebroid structures
on $A$ and $A^\ast$ and the 1-cocycles $\phi _0$ and $X_0$ satisfy
some compatibility conditions. When $\phi _0$ and $X_0$ are zero,
the definition reduces to that of a Lie bialgebroid. We also
remark that the theory of generalized Lie algebroids plays an
important role in the study of Lie brackets on affine bundles and
its application in the geometrical construction of Lagrangian-type
dynamics on affine bundles (see \cite{GGU,MMS}).

The aim of this paper is to integrate generalized Lie
bialgebroids, that is, to introduce the notion of a Jacobi
groupoid (a generalization of Poisson and contact groupoids), in
terms of groupoid morphisms, such that generalized Lie
bialgebroids to be considered as the infinitesimal invariants of
Jacobi groupoids.

As in the case of contact groupoids, we start with a Lie groupoid
$\gr$, a Jacobi structure $(\Jacobi )$ on $G$ and a multiplicative
function $\mult$. Then, as in the case of Poisson groupoids, we
consider the vector bundle morphism $\sostJ:T^\ast G\times\R\to
TG\times\R$ induced by the Jacobi structure $(\Jacobi )$. The
multiplicative function $\sigma$ induces, in a natural way, an
action of the tangent groupoid $TG\gpd TM$ over the canonical
projection $\pi _1:TM\times\R\to TM$ obtaining an action groupoid
$TG\times\R$ over $TM\times\R$. Thus, it is necessary to introduce
a suitable Lie groupoid structure in $T^\ast G\times\R$ over
$A^\ast G$ and this is the first important result of the paper. In
fact, we prove that:

$\bullet$ {\it If $AG$ is the Lie algebroid of an arbitrary Lie
groupoid $\gr$, $\mult$ is a multiplicative function,
$\bar{\pi}_G:T^\ast G\times \R\to G$ is the canonical projection
and $\eta _G$ is the canonical contact 1-form on $T^\ast
G\times\R$ then $(T^\ast G\times \R\gpd A^\ast G,\eta _G, \sigma
\circ \bar{\pi }_G)$ is a contact groupoid in such a way that the
Jacobi structure on $A^\ast G$ is just the linear Jacobi structure
$ (\Lambda_{(A^\ast G,\phi _0)},E_{(A^\ast G ,\phi_0)})$ induced
by the Lie algebroid $AG$ and the 1-cocycle $\phi _0$ which comes
from the multiplicative function $\sigma $ (see Theorems
\ref{contacto} and \ref{2.8}).}

Now, we will say that $(\gr ,\Jacobi ,\sigma )$ is a Jacobi
groupoid if the map $\sostJ :T^\ast G\times \R \to TG\times \R$ is
a Lie groupoid morphism over some map $\varphi _0:A^\ast G\to
TM\times \R$. Poisson and contact groupoids and other interesting
examples are Jacobi groupoids. In particular, Jacobi groupoids
$(\gr ,\Jacobi ,\sigma )$ where $M$ is a single point are just the
Lie groups studied in \cite{IM2}, whose infinitesimal invariants
are generalized Lie bialgebras.

On the other hand, if $(\gr ,\Jacobi ,\sigma )$ is a Jacobi
groupoid then we show that the vector bundle $A^\ast G$ admits a
Lie algebroid structure and the multiplicative function $\sigma$
(respectively, the vector field $E$) induces a 1-cocycle $\phi _0$
(respectively, $X_0$) on $AG$ (respectively, $A^\ast G$). Thus, a
first relation between Jacobi groupoids and generalized Lie
bialgebroids can be obtained and this is the second important
result of our paper:

$\bullet$ {\it If $(G\gpd M,\Jacobi ,\sigma )$ is a Jacobi
groupoid then $((AG,\phi _0),(A^\ast G,X_0))$ is a generalized Lie
bialgebroid (see Theorem \ref{bajada}).}

Finally, a converse of the above statement is the third important
result of the paper. More precisely, we prove the following
theorem:

$\bullet$ {\it Let $((AG,\phi _0),(A^\ast G,X_0))$ be a
generalized Lie bialgebroid where $AG$ is the Lie algebroid of an
$\s$-connected and $\s$-simply connected Lie groupoid $\gr$. Then,
there is a unique multiplicative function $\mult$ and a unique
Jacobi structure $(\Jacobi )$ on $G$ that makes
$(\gr,\Jacobi,\sigma)$ into a Jacobi groupoid with generalized Lie
bialgebroid $((AG,\phi _0),(A^\ast G,X_0))$ (see Theorem
\ref{subida}).}

The two previous results generalize those obtained by Mackenzie
and Xu \cite{MX,MX2} for Poisson groupoids and those obtained by
the authors \cite{IM2} for generalized Lie bialgebras.

The paper is organized as follows. In Section 2, we recall several
definitions and results about Jacobi manifolds, Lie algebroids and
Lie groupoids which will be used in the sequel. In Section 3, we
prove that a Lie groupoid $\gr$ (with Lie algebroid $AG$) and a
multiplicative function $\mult$ induce a Lie groupoid structure in
$TG\times\R$ over $TM\times\R$ and a contact groupoid structure in
$T^\ast G\times \R$ over $A^\ast G$. In Section 4, we introduce
the definition of a Jacobi groupoid, giving some examples, and we
prove some properties of these groupoids. In Section 5, we show
that generalized Lie bialgebroids are, in fact, the infinitesimal
invariants of Jacobi groupoids.

{\bf Notation:} If $M$ is a differentiable manifold, we will
denote by $C^\infty (M,\R)$ the algebra of $C^\infty$ real-valued
functions on $M$, by $\Omega ^k(M)$ the space of $k$-forms on $M$,
by $\mathfrak X (M)$ the Lie algebra of vector fields, by $\delta$
the de Rham differential on $\Omega ^\ast(M)=\oplus_{k}\Omega
^k(M)$, by ${\cal L}$ the Lie derivative operator and by $[\, ,\,
]$ the Schouten-Nijenhuis bracket (\cite{BV,V}). Moreover, if
$A\to M$ is a vector bundle over $M$ and $P \in \Gamma (\wedge ^2
A)$ is a section of $\wedge^2A\to M,$ we will denote by $\#_P
\colon A^\ast\to A$ the bundle map given by $\nu (\#_P(\omega
))=P(x) (\omega ,\nu )$, for $\omega, \nu \in A_x^\ast$, $A^\ast
_x$ being the fiber of $A^\ast$ over $x\in M$. We will also denote
by $\#_P\colon\Gamma (A^\ast)\to \Gamma (A)$ the corresponding
homomorphism of $C^\infty (M,\R)$-modules.
\section{Jacobi structures, Lie algebroids and Lie groupoids}
\subsection{Jacobi structures}
A {\em Jacobi structure} on a manifold $M$ is a pair $(\Lambda
,E)$, where $\Lambda$ is a 2-vector and $E$ is a vector field on
$M$ satisfying the following properties:
\begin{equation}\label{ecuaciones}
     [\Lambda ,\Lambda ]=2E\wedge \Lambda ,\hspace{1cm}
     [E,\Lambda ]=0.
\end{equation}
The manifold $M$ endowed with a Jacobi structure is called a {\em
Jacobi manifold}. A bracket of functions (the {\em Jacobi
bracket}) is defined by
   $$\{ f,g\} =\Lambda (\delta f,\delta g)+fE(g)-gE(f),$$
for all $f,g\in C^\infty (M,\R )$. In fact, the space $C^\infty
(M,\R)$ endowed with the Jacobi bracket is a {\em local Lie
algebra} in the sense of Kirillov (see \cite{K}). Conversely, a
structure of local Lie algebra on $C^\infty (M,\R)$ defines a
Jacobi structure on $M$ (see \cite{GL,K}). If the vector field $E$
identically vanishes then $(M,\Lambda )$ is a {\em Poisson
manifold} (see \cite{BV,Li1,V,We}).

Another interesting example of Jacobi manifolds comes from contact
manifolds. Let $M$ be a $2n+1$-dimensional manifold and $\eta$ a
1-form on $M$. We say that $(M,\eta )$ is a {\em contact manifold}
if $\eta \wedge (\delta \eta )^n\neq 0$ at every point (see, for
instance, \cite{LM,Li2}). A contact manifold $(M,\eta )$ is a
Jacobi manifold whose associated Jacobi structure $(\Lambda ,E)$
is given by $$\Lambda (\omega ,\nu )=\delta \eta (\flat _\eta
^{-1}(\omega ),\flat _\eta ^{-1}(\nu )), \qquad E =\flat _\eta
^{-1}(\eta ),$$ for $\omega ,\nu \in \Omega ^1(M)$, $\flat _\eta
:\mathfrak X (M) \to \Omega ^1(M)$ being the isomorphism of
$C^\infty (M,\R)$-modules defined by $\flat _\eta (X)=i(X)\delta
\eta +\eta (X)\eta$. Note that $E$ is the {\em Reeb vector field}
of $M$ which is characterized by the conditions $i(E)\eta =1$ and
$i(E)\delta \eta =0$. Moreover, $$\flat _\eta \inv (\omega
)=-\sostP (\omega )+\omega (E)E,\mbox{ for }\omega \in \Omega
^1(M).$$ Jacobi manifolds were introduced by Lichnerowicz
\cite{Li2} (see also \cite{DLM,GL}).
\begin{remark}\label{poisonizacion}
{\rm Let $(\Jacobi )$ be a 2-vector and a vector field on a
manifold $M$. Then, we can consider the 2-vector $\tilde{\Lambda}$
given by
\begin{equation}\label{homog}
\tLambda=e^{-t}\Big ( \Lambda +\frac{\partial }{\partial t}\wedge
E \Big ) ,
\end{equation}
where $t$ is the usual coordinate on $\R$. The 2-vector $\tLambda$
is homogeneous with respect to the vector field $\fract$, that is,
${\cal L}_{\fract }\tLambda =-\tLambda$. In fact, if $\tLambda$ is
a 2-vector on $M\times\R$ such that ${\cal L}_{\fract }\tLambda
=-\tLambda$ then there exists a 2-vector $\Lambda$ and a vector
field $E$ on $M$ such that $\tLambda$ is given by (\ref{homog}).
Moreover, $(\Jacobi )$ is a Jacobi structure on $M$ if and only if
$\tLambda$ defines a Poisson structure on $M\times \R$ (see
\cite{Li2}). The manifold $M\times \R$ endowed with the structure
$\tilde{\Lambda}$ is called the {\em Poissonization of the Jacobi
manifold} $(M,\Lambda ,E)$. If $(\Jacobi )$ is a Jacobi structure
on $M$ induced by a contact 1-form $\eta$ then the corresponding
Poisson structure $\tilde{\Lambda}$ on $M\times\R$ is
non-degenerate and is associated with the symplectic 2-form
$\tilde{\Omega}=e^{t}\Big ( \delta \eta +\delta t\wedge \eta \Big
)$. }
\end{remark}
Before finishing this Section, we will give a definition which
will be useful in the following.
\begin{definition}
Let $S$ be a submanifold of a manifold $M$ and $\Lambda$ be an
arbitrary 2-vector. $S$ is said to be {\em coisotropic} (with
respect to $\Lambda$) if $\sostP( (T_xS)^\circ )\subseteq T_xS$,
for $x\in S$, $(T _xS)^\circ$ being the annihilator space of
$T_xS$.
\end{definition}
\begin{remark}
{\rm If $\Lambda$ (respectively, $(\Jacobi )$) is a Poisson
structure (respectively, a Jacobi structure) on $M$ then we
recover the notion of a coisotropic submanifold of the Poisson
manifold $(M,\Lambda )$ \cite{LM,We2} (respectively, coisotropic
submanifold of a Jacobi manifold $(M,\Jacobi )$ \cite{ILMM}) .}
\end{remark}
\subsection{Lie algebroids}\label{algebroides}
A {\em Lie algebroid} $A$ over a manifold $M$ is a vector bundle
$A$ over $M$ together with a Lie bracket $\lcf\, ,\,\rcf$ on the
space $\Gamma (A)$ of the global cross sections of $A\to M$ and a
bundle map $\rho \colon A \to TM$, called the {\em anchor map},
such that if we also denote by $\rho \colon \Gamma (A) \to
\mathfrak X (M)$ the homomorphism of $C^\infty (M,\R )$-modules
induced by the anchor map then:
\begin{enumerate}
\item[(i)] $\rho \colon (\Gamma (A),\lcf \, ,\, \rcf )\to (\mathfrak X
(M),[\, ,\, ])$ is a Lie algebra homomorphism and
\item[(ii)] for all $f\in C^\infty (M,\R)$ and for all $X ,Y \in
\Gamma (A)$, one has $$ \lcf X ,fY \rcf=f\lcf X,Y\rcf +(\rho  (X
)(f))Y. $$
\end{enumerate}
The triple $\alg$ is called a Lie algebroid over $M$ (see
\cite{Mk,Pr}).

A real Lie algebra of finite dimension is a Lie algebroid over a
point. Another example of a Lie algebroid is the triple $(TM,[\,
,\, ],Id)$, where $M$ is a differentiable manifold and $Id\colon
TM\to TM$ is the identity map.

If $A$ is a Lie algebroid, the Lie bracket on $\Gamma (A)$ can be
extended to the so-called {\em Schouten bracket} $\lcf \, ,\,
\rcf$ on the space $\Gamma (\wedge ^\ast A)=\oplus_k\Gamma (\wedge
^kA)$ of multi-sections of $A$ in such a way that $(\oplus_k\Gamma
(\wedge ^kA),\wedge ,\lcf \, ,\, \rcf )$ is a graded Lie algebra.
In fact, the Schouten bracket satisfies the following properties
$$\lcf X ,f\rcf =\rho (X )(f),\quad \lcf P,Q \rcf =(-1)^{pq}\lcf
Q,P\rcf ,$$ $$\lcf P,Q\wedge R\rcf =\lcf P,Q\rcf \wedge
R+(-1)^{q(p+1)}Q\wedge \lcf P,R\rcf ,$$ $$ (-1)^{pr}\lcf \lcf P,Q
\rcf ,R\rcf +
   (-1)^{qr}\lcf \lcf R,P \rcf ,Q\rcf +
   (-1)^{pq}\lcf \lcf Q,R \rcf ,P\rcf =0,$$
for $X \in \Gamma(A)$, $f\in C^\infty (M,\R)$, $P\in \Gamma
(\wedge ^pA)$, $Q\in \Gamma (\wedge ^{q}A)$ and  $R\in \Gamma
(\wedge ^{r}A)$ (see \cite{V}).
\begin{remark}\label{signo}
{\rm The definition of Schouten bracket considered here is the one
given in \cite{V} (see also \cite{BV,Li1}). Some authors, see for
example \cite{K-S}, define the Schouten bracket in another way. In
fact, the relation between the Schouten bracket $\lcf\;,\;\rcf'$
in the sense of \cite{K-S} and the Schouten bracket
$\lcf\;,\;\rcf$ in the sense of \cite{V} is the following one. If
$P\in \Gamma(\wedge^{p}A)$ and $Q\in \Gamma(\wedge^\ast A)$, then
$\lcf P,Q\rcf'=(-1)^{p+1}\lcf P,Q\rcf$. }
\end{remark}
On the other hand, imitating the de Rham differential on the space
$\Omega ^\ast (M)$, we define the {\em differential of the Lie
algebroid} $A$, $d\colon \Gamma (\wedge ^kA^\ast )$ $\to \Gamma
(\wedge ^{k+1}A^\ast )$, as follows. For $\omega\in \Gamma (\wedge
^kA^\ast )$ and $X_0,\ldots ,X_k\in \Gamma (A)$,
\begin{equation}\label{diferencial}
\begin{array}{ccl}
d\omega (X_0,\ldots,X_k)&=&\kern-5pt\displaystyle
\sum_{i=0}^k(-1)^i\rho(X_i)(\omega (X_0,\ldots ,\hat{X}_i,\ldots
,X_k))\\ &+&\kern-5pt\displaystyle\sum_{i<j}(-1)^{i+j}\omega (\lcf
X_i,X_j\rcf ,X_0,\ldots ,\hat{X}_i,\ldots ,\hat{X}_j,\ldots ,X_k).
\end{array}
\end{equation} Moreover, since $d^2=0$, we have the corresponding
cohomology spaces. This cohomology is the {\em Lie algebroid
cohomology with trivial coefficients} (see \cite{Mk}).

Using the above definitions, it follows that a 1-cochain $\phi \in
\Gamma (A^\ast)$ is a 1-cocycle if and only if $$ \phi \lcf X ,Y
\rcf =\rho (X )(\phi (Y))-\rho (Y)(\phi (X )), $$ for all $X ,Y
\in\Gamma (A)$.

Next, we will consider some examples of Lie algebroids which will
be important in the following.

{\bf 1.-} {\em The Lie algebroid} $(TM\times \R ,\makebox{{\bf
[}\, ,\, {\bf ]}}, \pi)$

If $M$ is a differentiable manifold, then the triple $(TM\times \R
,\makebox{{\bf [}\, ,\, {\bf ]}}, \pi)$ is a Lie algebroid over
$M$, where $\pi \colon TM\times \R \to TM$ is the  canonical
projection over the first factor and $\makebox{{\bf [}\, ,\, {\bf
]}}$ is the bracket given by (see \cite{Mk,NVQ})
\begin{equation}\label{corchetedeprimerorden}
\makebox{{\bf [}} (X,f),(Y,g)\makebox{{\bf ]}}=([X,Y],X(g)-Y(f)),
\end{equation}
for $(X,f),(Y,g )\in \mathfrak X (M)\times C^\infty (M,\R )\cong
\Gamma (TM\times \R)$.

{\bf 2.-} {\em The Lie algebroid $(T^\ast M \kern-1pt \times \R
,\lcf ,\rcf _{(\Lambda ,E)},\kern-2pt \widetilde{\#}_{(\Lambda
,E)})$ associated with a Jacobi manifold \kern-1pt $(\kern-2pt
M,\Lambda ,E)$}

A Jacobi manifold $(M,\Lambda ,E)$ has an associated Lie algebroid
$(T^\ast M \times \R ,\lcf ,\rcf _{(\Lambda ,E)},$
$\widetilde{\#}_{(\Lambda ,E)})$, where $\lcf \, ,\, \rcf
_{(\Lambda ,E)}$ and $\widetilde{\#}_{(\Lambda,E)}$ are defined by
\begin{equation}\label{ecjacobi}
\begin{array}{cll}
\kern-14pt\lcf (\omega ,f),(\nu ,g)\rcf _{(\Lambda
,E)}&\kern-11pt=&\kern-11pt({\cal L}_{\#_{\Lambda}(\omega
)}\nu\kern-2pt-\kern-2pt{\cal L}_{\#_{\Lambda}(\nu
)}\omega\kern-2pt -\kern-2pt \delta (\Lambda (\omega ,\nu
))\kern-3pt+\kern-3ptf{\cal L}_{E}\nu\kern-2pt-\kern-2ptg {\cal
L}_{E}\omega
\\&\kern-99pt &\kern-99pt -i(E)(\omega \kern-1pt\wedge\kern-1pt \nu),
\Lambda (\nu ,\omega )\kern-3pt+\kern-3pt\#_{\Lambda}(\omega
)(g)\kern-3pt-\kern-3pt\#_{\Lambda}(\nu
)(f)\kern-3pt+\kern-3ptfE(g)\kern-3pt-\kern-3ptg E(f)),\\[15pt]
\kern-17pt\widetilde{\#}_{(\Lambda ,E)}(\omega
,f)&\kern-19pt=&\kern-10pt\#_{\Lambda}(\omega )+fE,
\end{array}
\end{equation}
for $(\omega ,f),(\nu ,g)\in \Omega ^1(M)\times C^\infty (M,\R
)\cong \Gamma (T^\ast M\times \R)$, ${\cal L}$ being the Lie
derivative operator (see \cite{KS}). In the particular case when
$(M,\Lambda )$ is a Poisson manifold we recover, by projection,
the Lie algebroid $(T^\ast M,\lcf \, ,\, \rcf _{\Lambda},
\#_{\Lambda})$, where $\lcf \, ,\, \rcf_{\Lambda}$ is the bracket
of 1-forms defined by $\lcf \omega ,\nu \rcf _{\Lambda}={\cal
L}_{\#_{\Lambda}(\omega )}\nu-{\cal L}_{\#_{\Lambda}(\nu )}\omega
-\delta (\Lambda (\omega ,\nu )),$ for $\omega ,\nu \in \Omega
^1(M)$ (see \cite{BV,CDW,F,V}).

{\bf 3.-} {\em Action of a Lie algebroid on a smooth map}

Let $\alg$ be a Lie algebroid over a manifold $M$ and $\pi :P\to
M$ be a smooth map. An action of $A$ on $\pi :P\to M$ is a
$\R$-linear map $$\ast :\Gamma (A)\to \mathfrak X (P),\quad X\in
\Gamma (A)\mapsto X^\ast \in \mathfrak X (P),$$ such that:
$$(fX)^\ast =(f\circ \pi)X^\ast,\quad \lcf X, Y\rcf ^\ast =[
X^\ast ,Y^\ast ],\quad \pi ^p_\ast(X^\ast (p))=\rho (X (\pi
(p))),$$ for $f\in C^\infty (M,\R)$, $X,Y\in \Gamma (A)$ and $p\in
M$. If $\ast :\Gamma (A)\to \mathfrak X (P)$ is an action of $A$
on $\pi :P\to M$ and $\tau :A\to M$ is the bundle projection then
the pullback vector bundle of $A$ over $\pi$ $$\pi ^\ast A=\{
(a,p)\in A\times P \,/\, \tau (a)=\pi (p)\}$$ is a Lie algebroid
over $P$ with the Lie algebroid structure $(\lcf \, ,\, \rcf
_{\pi},\rho _{\pi })$ which is characterized by $$\rho _\pi
(X)(p)=X^\ast (p),\quad \lcf X,Y\rcf _\pi = \lcf X,Y\rcf \circ \pi
,$$ for $X,Y\in \Gamma (A)$ and $p\in P$. The triple $(\pi ^\ast
A,\lcf \, ,\, \rcf _{\pi},\rho _{\pi })$ is called the {\em action
Lie algebroid of $A$ on $\pi$} and it is denoted by $A\ltimes \pi$
or $A\ltimes P$ (see \cite{HM}).

{\bf 4.-} {\em The Lie algebroid associated with a linear Poisson
structure}

Let $\tau :A\to M$ be a vector bundle on a manifold $M$. Then, it
is clear that there exists a bijection between the space $\Gamma
(A^\ast )$ of the sections of the dual bundle $\tau ^\ast :A^\ast
\to M$ and the set ${\cal L}(A)$ of real functions on $A$ which
are linear on each fiber $$\Gamma (A^\ast )\to {\cal L}(A),\quad
\nu \to \tilde{\nu }.$$ Now, suppose that $\Lambda$ is a linear
Poisson structure on $A$ with Poisson bracket $\{ \, ,\, \}$. This
means that the Poisson bracket of two linear functions on $A$ is
again a linear function. This fact implies that the Poisson
bracket of a linear function on $A$ and a basic function is a
basic function. Moreover, one may define a Lie algebroid structure
$(\lcf \, ,\, \rcf ,\rho )$ on $\tau ^\ast :A^\ast \to M$ which is
characterized by
\begin{equation}\label{P-lineal}
\widetilde{\lcf \nu ,\mu \rcf}=\{ \tilde{\nu },\tilde{\mu}\},
\quad \rho (\nu )(f)\circ \tau  = \{ \tilde{\nu },f\circ \tau \},
\end{equation}
for $\nu ,\mu \in \Gamma (A^\ast )$ and $f\in C^\infty (M,\R )$
(see \cite{CDW,Co}). Conversely, if $A$ is a vector bundle over
$M$ and the dual bundle $A^ \ast$ admits a Lie algebroid structure
$\estalg$ then one may define a linear Poisson bracket $\{ \,
,\,\}$ on $A$ in such a way that (\ref{P-lineal}) holds.

{\bf 5.-} {\em The tangent Lie algebroid}

Let $(M,\Lambda )$ be a Poisson manifold. Then, the complete lift
$\Lambda ^c$ of $\Lambda$ to the tangent bundle $TM$ defines a
linear Poisson structure on $TM$ (see \cite{Co,S}). $\Lambda ^c$
is called the {\em tangent Poisson structure}.

Now, suppose that $\tau :A\to M$ is a Lie algebroid over a
manifold $M$ and that $p :A^\ast \times _M A\to \R$ is the natural
pairing. Then, $TA$ and $TA^\ast$ are vector bundles over $TM$ and
$p$ induces a non-degenerate pairing $TA^\ast \times _{TM} TA\to
\R$. Thus, we get an isomorphism between the vector bundles $TA$
and $(TA^\ast )^\ast$. Therefore, the dual bundle to $TA\to TM$
may be identified with $TA^\ast \to TM$. On the other hand, since
$A^\ast$ is a Poisson manifold, we have that $TA^\ast$ admits a
linear Poisson structure. Consequently, the vector bundle $TA\to
TM$ is a Lie algebroid which is called the {\em tangent Lie
algebroid} to $A$ (for more details, see \cite{Co2,MX}).
\subsection{Lie groupoids}\label{grupoides}
A {\em groupoid} consists of two sets $G$ and $M$, called
respectively the {\em groupoid} and the {\em base}, together with
two maps $\s$ and $\t$ from $G$ to $M$, called respectively the
{\em source} and {\em target} projections, a map $\epsilon :M\to
G$, called the {\em inclusion}, a partial multiplication $m :
G^{(2)}=\{ (g,h)\in G\times G / \s (g)=\t (h)\} \to G$ and a map
$\iota: G\to G$, called the {\em inversion}, satisfying the
following conditions:
\begin{enumerate}
\item $\s (m(g,h))=\s (h)$ and $\t (m(g,h))=\t (g)$, for all
$(g,h)\in G^ {(2)}$,
\item $m(g,m(h,k))=m(m(g,h),k)$, for all $g,h,k\in G$
such that $\s (g)=\t (h)$ and $\s(h) =\t (k)$,
\item $\s (\epsilon (x)) =x$ and $\t (\epsilon (x))=x$, for all
$x\in M$,
\item $m(g,\epsilon (\s (g)))=g$ and $m(\epsilon (\t
(g)),g)=g$, for all $g\in G$,
\item $m(g,\iota (g))=\epsilon (\t (g))$ and $m(\iota (g),g)
=\epsilon (\s (g))$, for all $g\in G$.
\end{enumerate}
A groupoid $G$ over a base $M$ will be denoted by $\gr$. Given two
groupoids $G_1\gpd M_1$ and $G_2\gpd M_2$, a {\em morphism of
groupoids} is a pair of maps $\Phi :G_1\to G_2$ and $\Phi
_0:M_1\to M_2$ which commute with all the structural functions of
$G_1$ and $G_2$, i.e., $\s _2 \circ \Phi =\Phi _0\circ \s_1$, $\t
_2 \circ \Phi =\Phi _0\circ \t_1$ and $\Phi (g_1h_1)=\Phi
(g_1)\Phi (h_1)$, for $(g_1,h_1)\in G_1^{(2)}$ (for more details,
see \cite{Mk}).

If $G$ and $M$ are manifolds, $\gr$ is a {\em Lie groupoid} if:
\begin{enumerate}
\item $\s$ and $\t$ are differentiable submersions.
\item $m$, $\epsilon$ and $\iota$ are differentiable maps.
\end{enumerate}
From now on, we will usually write $gh$ for $m(g,h )$, $g\inv$ for
$\iota (g)$ and $\tilde{x}$ for $\epsilon (x)$. Moreover, if $x\in
M$ then $G _x=\s\inv (x)$ (respectively, $G ^x=\t\inv (x)$) will
be said the $\s$-fiber (resp., the $\t$-fiber) of $x$.
Furthermore, since $\epsilon$ is an inmersion, we will identify
$M$ with $\epsilon (M)$.

Next, we will recall some notions related with Lie groupoids which
will be useful in the following (for more details, see \cite{Mk}).
\begin{definition}
Let $\gr$ be a Lie groupoid over a manifold $M$. For $U\subseteq
M$ open, a {\em local bisection} (or {\em local admissible
section}) of $G$ on $U$ is a smooth map ${\cal K} :U\to G$ which
is right-inverse to $\t$ and for which $\s \circ {\cal K}:U\to \s
({\cal K}(U))$ is a diffeomorphism from $U$ to the open set $\s
({\cal K}(U))$ in $M$. If $U=M$, ${\cal K}$ is a global bisection
or simply a bisection.
\end{definition}
The existence of local bisections through any point $g\in G$ is
always guaranteed.

If ${\cal K} :U\to G$ is a local bisection with $V=(\s \circ {\cal
K}) (U)$, the local {\em left-translation} and {\em
right-translation} induced by ${\cal K}$ are the maps $L_{\cal
K}:\t \inv (V)\to \t \inv (U)$ and $R_{\cal K}:\s \inv (U)\to \s
\inv (V)$, defined by
\[
L_{\cal K}(g)={\cal K}((\s \circ {\cal K})\inv (\t (g)))\, g,\quad
R_{\cal K} (h)=h\, {\cal K}(\s (h)),
\]
for $g\in \t \inv (V)$ and $h\in \s \inv (U)$.
\begin{remark}
{\rm If $y_0\in U$ and ${\cal K} (y_0)=g_0$, $\s (g_0)=x_0$ then
the restriction of $L_{\cal K}$ to $G^{x_0}$ is the
left-translation by $g_0$ $$L_{g_0}:G^{x_0}\to G^{y_0}, \quad
h\mapsto L_{g_0}(h)=g_0h.$$ In a similar way, the restriction of
$R_{\cal K}$ to $G_{y_0}$ is the right-translation by $g_0$
$$R_{g_0}:G_{y_0}\to G_{x_0}, \quad g\mapsto R_{g_0}(g)=gg_0.$$}
\end{remark}
A multivector field $P$ on $G$ is said to be {\em left-invariant}
(respectively, {\em right-invariant}) if it is tangent to the
fibers of $\t$ (respectively, $\s$) and $P(gh)=(L_{\cal K})_\ast
^h (P(h))$ (respectively, $P(gh)=(R_{\cal K})_\ast ^g (P(g))$) for
$g,h\in G$ and ${\cal K}:U\to G$ any local bisection through $h$
(res\-pectively, $g$). If $P$ and $Q$ are two left-invariant
(respectively, right-invariant) multivector fields on $G$ then
$[P,Q]$ is again left-invariant (respectively, right-invariant).

Now, we will recall the definition of the Lie algebroid associated
with a Lie groupoid.

Suppose that $\gr$ is a Lie groupoid. Then, we may consider the
vector bundle $AG\to M$, whose fiber at a point $x\in M$ is
$A_xG=T_{\tilde{x}}G^x$. It is easy to prove that there exists a
bijection between the space $\Gamma (AG)$ and the set of
left-invariant (respectively, right-invariant) vector fields on
$G$. If $X$ is a section of $AG$, the corresponding left-invariant
(respectively, right-invariant) vector field on $G$ will be
denoted by $\izq{X}$ (respectively, $\der{X}$). Using the above
facts, we may introduce a Lie algebroid structure $(\lcf \, ,\,
\rcf ,\rho )$ on $AG$, which is defined by, for $X,Y\in \Gamma
(AG)$ and $x\in M$,
\begin{equation}\label{inv-izq}
\izq{\lcf X, Y\rcf}=[\izq{X},\izq{Y}],\quad \rho (X)(x)=\s _\ast
^{\tilde{x}}(X(x)).
\end{equation}

\vspace{-.5cm}
\begin{remark}
{\rm There exists a bijection between the space $\Gamma (\wedge
^k(AG))$ and the set of left-invariant (respectively,
right-invariant) $k$-vector fields. If $P$ is a section of $\wedge
^k (AG)$, we will denote by $\izq{P}$ (respectively, $\der{P}$)
the corresponding left-invariant (respectively, right-invariant)
$k$-vector field on $G$. Moreover, if $P,Q\in \Gamma (\wedge ^\ast
(AG))$, we have that
\begin{equation}\label{multi-inv-izq}
\izq{\lcf P, Q\rcf}=[\izq{P},\izq{Q}].
\end{equation}
}
\end{remark}

\begin{examples}\label{ej-grupoides}
{\rm {\bf 1.-} {\em Lie groups}

Any Lie group $G$ is a Lie groupoid over $\{\mathfrak e\}$, the
identity element of $G$. The Lie algebroid associated with $G$ is
just the Lie algebra $\mathfrak g$ of $G$.

{\bf 2.-} {\em The banal groupoid}

Let $M$ be a differentiable manifold. The product manifold
$M\times M$ is a Lie groupoid over $M$ in the following way: $\s$
is the projection onto the second factor and $\t$ is the
projection onto the first factor; $\epsilon (x)=(x,x)$ for all
$x\in M$ and $m((x,y),(y,z))=(x,z)$. $M\times M \gpd M$ is called
the {\em banal groupoid}. The Lie algebroid associated with the
banal groupoid is the tangent bundle $TM$ of $M$.

{\bf 3.-} {\em The direct product of Lie groupoids}

If $G_1\gpd M_1$ and $G_2\gpd M_2$ are Lie groupoids, then
$G_1\times G_2\gpd M_1\times M_2$ is a Lie groupoid in a natural
way.

{\bf 4.-} {\em Action groupoids}

Let $\gr$ be a Lie groupoid and $\pi :P\to M$ be a smooth map. If
$P\ast G =\{ (p,g)\in P\times G\, /\, \pi (p)=\t (g) \}$ then a
right action of $G$ on $\pi$ is a smooth map $$P\ast G \to P,\quad
(p,g)\mapsto p\cdot g,$$ which satisfies the following relations
$$\begin{array}{l} \pi (p\cdot g)=\s (g), \mbox{ for all }
(p,g)\in P\ast G,\\ (p \cdot g)\cdot h=p\cdot (gh), \mbox { for
all } (g,h)\in G^{(2)} \mbox{ and } (p,g) \in P\ast G,\\ p\cdot
\widetilde{\pi (p)}=p, \mbox{ for all }p\in P.
\end{array}$$
Given such an action one constructs the {\em action groupoid}
$P\ast G\gpd P$ by defining $$\s '(p,g)=p\cdot g, \quad \t
'(p,g)=p,$$ $$m'((p,g),(q,h))=(p, gh),\mbox{ if }\, q=p\cdot g,$$
$$\epsilon '(p)=(p,\epsilon (\pi (p))),\quad \iota '(p,g)=(p\cdot
g,g\inv ).$$ Now, if $p\in P$, we consider the map $\pi _p:G^{\pi
(p)}\to P$ given by $$\pi _p(g)=p\cdot g.$$ Then, if $AG$ is the
Lie algebroid of $G$, the $\R$-linear map $$\ast :\Gamma (AG)\to
\mathfrak X(P),\;\;\;\; X\in \Gamma (AG)\mapsto X^\ast \in
\mathfrak X(P),$$ defined by
\begin{equation}\label{1.9I}
X^\ast (p)=(\pi _p)_\ast ^{\widetilde{\pi (p)}}(X(\pi (p))),\mbox{
for all }p\in P,
\end{equation}
induces an action of $AG$ on $\pi :P\to M$. In addition, the Lie
algebroid associated with the Lie groupoid $P\ast G\gpd P$ is the
action Lie algebroid $AG\ltimes \pi$ (for more details, see
\cite{HM}).

{\bf 5.-} {\em The tangent groupoid}

Let $\gr$ be a Lie groupoid. Then, the tangent bundle $TG$ is a
Lie groupoid over $TM$. The projections $\s ^T$, $\t ^T$, the
partial multiplication $\oplus _{TG}$, the inclusion $\epsilon ^T$
and the inversion $\iota ^T$ are defined by
\begin{equation}\label{TG}
\begin{array}{l}
\s ^T(X_g)=\s _\ast ^g(X_g),\quad \t ^T(X_g)=\t _\ast ^g(X_g),
\mbox{ for }X_g\in T_gG,\\ X_g\oplus _{TG}Y_h= m_\ast
^{(g,h)}(X_g,Y_h), \mbox{ for } (X_g,Y_h)\in
(TG)_{(g,h)}^{(2)}=T_{(g,h)}G^{(2)},\\ \epsilon ^T(X_x)=\epsilon
_\ast ^x(X_x),\mbox{ for }X_x\in T_xM,\\ \iota ^T (X_g)=\iota
_\ast ^g (X_g),\mbox{ for }X_g\in T_gG.\end{array}
\end{equation}
In \cite{Xu} it has been given an explicit expression for the
multiplication $\oplus _{TG}$. If $\s _\ast ^g(X_g)=\t _\ast
^h(X_h)=W_x$, $x=\s (g)=\t (h)$, then
\begin{equation}\label{sumaTG}
X_g\oplus _{TG} Y_h = (L_{{\cal X}})_\ast ^h(Y_h) + (R_{{\cal
Y}})_\ast ^g(X_g)- (L_{{\cal X}})_\ast ^h((R_{{\cal Y}})_\ast
^{\tilde{x}}(\epsilon _\ast ^x (W))),
\end{equation}
where ${\cal X}, {\cal Y}$ are any (local) bisections of $G$ with
${\cal X}(x) = g$ and ${\cal Y}(x) = h$. The tangent Lie algebroid
$TAG\to TM$ is just the Lie algebroid associated with the tangent
groupoid $TG\gpd TM$ (for more details, see \cite{MX}).
\begin{remark}
{\rm If $G$ is a Lie group then, from (\ref{sumaTG}), it follows
that
\begin{equation}\label{mult'}
X_g\oplus _{TG} Y_h = (L_g)_\ast ^h(Y_h) + (R_h)_\ast
^g(X_g),\mbox{ for }X_g\in T_gG \mbox{ and }Y_h\in T_hG.
\end{equation}}
\end{remark}
{\bf 6.-} {\em The cotangent groupoid}

Let $\gr$ be a Lie groupoid. If $A^\ast G$ is the dual bundle to
$AG$ then the cotangent bundle $T^\ast G$ is a Lie groupoid over
$A^\ast G$. The projections $\tilde{\s}$ and $\tilde{\t}$, the
partial multiplication $\oplus _{T^\ast G}$, the inclusion
$\tilde{\epsilon}$ and the inversion $\tilde{\iota}$ are defined
as follows,
\begin{equation}\label{T*G}
\begin{array}{l} \tilde{\s}(\omega _g)(X)=\omega _g
((L_g)_\ast ^{\widetilde{\s (g)}} (X)),\mbox{ for } \omega _g\in
T^\ast _gG \mbox{ and } X\in A_{\s (g)}G,\\ \tilde{\t}(\nu
_h)(Y)=\nu _h ((R_h)_\ast ^{\widetilde{\t (h)}}(Y-\epsilon _\ast
^{\t (h)} (\s _\ast ^{\widetilde{\t (h)}}(Y)))),\mbox{ for } \nu
_h\in T^\ast _hG\mbox{ and }Y\in A_{\t(h)}G,\\ (\omega _g\oplus
_{T^\ast G}\nu _h)(X_g\oplus _{TG} Y_h)=\omega _g(X_g)+\nu _h
(Y_h),\mbox{ for }(X_g,Y_h)\in T_{(g,h)}G^{(2)},\\
\tilde{\epsilon}(\omega _x)(X_{\tilde{x}})=\omega
_x(X_{\tilde{x}}-\epsilon _\ast ^x (\t ^{\tilde{x}}_\ast
(X_{\tilde{x}}))), \mbox{ for }\omega _x\in A^\ast _xG,
X_{\tilde{x}}\in T_{\tilde{x}}G \mbox{ and }x\in M,\\
\tilde{\iota} (\omega _g)(X_{g\inv})=-\omega _g(\iota  _\ast
^{g\inv} (X_{g\inv})), \mbox{ for } \omega _g\in T^\ast _gG \mbox{
and } X_{g\inv }\in T_{g\inv}G.\end{array}
\end{equation}
Note that $\tilde{\epsilon}( A^\ast G )$ is just the conormal
bundle of $M\cong \epsilon (M)$ as a submanifold of $G$.

On the other hand, since $A^\ast G$ is a Poisson manifold, the
cotangent bundle $T^\ast (A^\ast G)$ is a Lie algebroid. In fact,
the Lie algebroid of the cotangent Lie groupoid $T^\ast G\gpd
A^\ast G$ may be identified with $T^\ast (A^\ast G)$ (for more
details, see \cite{CDW,MX}).
\begin{remark}
{\rm If $G$ is a Lie group and $\omega _g\in T^\ast _gG$, $\nu
_h\in T^\ast _hG$ satisfy $\tilde{\s}(\omega _g)=\tilde{\t}(\nu
_h)$ then, from (\ref{mult'}), it follows that
\begin{equation}\label{mult*}
\omega _g\oplus _{T^\ast G}\nu _h =\displaystyle \frac{1}{2}\Big
\{ ((R_{h\inv})_\ast ^{gh})^\ast (\omega _g)+((L_{g\inv})_\ast
^{gh})^\ast (\nu _h)\Big \}.
\end{equation}
}
\end{remark}
}
\end{examples}
\subsection{Generalized Lie bialgebroids}\label{glbi}
In this Section, we will recall the definition of a generalized
Lie bialgebroid. First, we will exhibit some results about the
differential calculus on Lie algebroids in the presence of a
1-cocycle (for more details, see \cite{IM}).

If $\alg$ is a Lie algebroid over $M$ and, in addition, we have a
1-cocycle $\phi _0\in \Gamma (A ^\ast)$ then the usual
representation of the Lie algebra $\Gamma (A)$ on the space
$C^\infty (M,\R )$ can be modified and a new representation is
obtained. This representation is given by $\rho _{\phi
_0}(X)(f)=\rho (X)(f)+\phi _0 (X)f$, for $X\in \Gamma (A)$ and
$f\in C^\infty (M,\R )$. The resulting cohomology operator
$d_{\phi _0}$ is called the $\phi _0$-differential of $A$ and its
expression, in terms of the differential $d$ of $A$, is
\begin{equation}\label{ndif}
d_{\phi _0}\omega =d\omega +\phi _0\wedge \omega
\end{equation}
for $\omega \in \Gamma (\wedge ^k A^\ast )$. The $\phi
_0$-differential of $A$ allows us to define, in a natural way, the
$\phi _0$-Lie derivative by a section $X\in \Gamma (A)$, $({\cal
L}_{\phi _0})_X\colon \Gamma (\wedge ^kA^\ast)\to \Gamma (\wedge
^kA^\ast)$, as the commutator of $d_{\phi _0}$ and the contraction
by $X$, that is, $({\cal L}_{\phi _0})_X  =d_{\phi _0}\circ
i(X)+i(X)\circ d_{\phi _0}$ (for the general definition of the
differential and the Lie derivative associated with a
representation of a Lie algebroid on a vector bundle, see
\cite{Mk}).

On the other hand, imitating the definition of the Schouten
bracket of two multilinear first-order differential operators on
the space of $C^\infty$ real-valued functions on a manifold $N$
(see \cite{BV}), we introduced the $\phi _0$-Schouten bracket of a
$p$-section $P$ and a $p'$-section $P'$ as the ($p+p'-1$)-section
given by
\begin{equation}\label{Schouext}
\lcf P, P' \rcf_{\phi _0}=\lcf P, P' \rcf + (-1)^{p+1}(p-1)P\wedge
(i(\phi _0)P')-(p'-1)(i(\phi_0)P)\wedge P',
\end{equation}
where $\lcf \, ,\, \rcf$ is the usual Schouten bracket of $A$
(some properties of the $\phi _0$-Schouten bracket were obtained
in \cite{IM}). Moreover, using the $\phi _0$-Schouten bracket, we
can define the $\phi _0$-Lie derivative of $P\in \Gamma (\wedge ^k
A)$ by $X \in \Gamma (A)$ as
\begin{equation}\label{Lieext}
({\cal L}_{\phi _0})_X (P)=\lcf X ,P\rcf _{\phi _0}.
\end{equation}

\vspace{-.6cm}
\begin{remark}
{\rm The product manifold $\bar{A}=A\times T\R$ is a vector bundle
over $M\times \R$ and one may define a Lie algebroid structure
$(\lcf \, ,\, \rcf \, \bar{ },\bar{\rho })$ on $\bar{A}$, where
$\lcf \, ,\, \rcf \, \bar{ }$ is the obvious product Lie bracket
and $\bar{\rho}=\rho \times id: \bar{A}\to TM\times T\R$. The
direct sum $\Gamma (\wedge ^pA)\oplus \Gamma (\wedge ^{p-1} A)$ is
a subspace of $\Gamma (\wedge ^p\bar{A})$ and we may consider the
monomorphism of $C^\infty (M,\R )$-modules $\bar{U}_{\phi
_0}:\Gamma (\wedge ^p A)\to\Gamma (\wedge ^p \bar{A})$ given by
$\bar{U}_{\phi _0}(P)=(e^{-(p-1)t}P, e^{-(p-1)t} i(\phi _0)(P))$.
Then, it is easy to prove that $\bar{U}_{\phi _0}(\lcf P,P'\rcf
_{\phi _0})=\lcf \bar{U}_{\phi _0}(P),\bar{U}_{\phi _0}(P')\rcf
\,\bar{ }$, for $P\in \Gamma (\wedge ^p A)$ and $P'\in \Gamma
(\wedge ^{p'}A)$ (see \cite{GM}). }
\end{remark}
Now, suppose that $\alg$ is a Lie algebroid and that $\phi _0\in
\Gamma (A^\ast )$ is a 1-cocycle. Assume also that the dual bundle
$A^\ast$ admits a Lie algebroid structure $(\lcf \, ,\, \rcf _\ast
,\rho _\ast )$ and that $X_0\in \Gamma (A)$ is a 1-cocycle. The
pair $((A,\phi _0),(A^\ast ,X_0))$ is a {\em generalized Lie
bialgebroid} if
\begin{equation}\label{condcomp}
\begin{array}{l}
 d_\ast{}_{X_0}\lcf X ,Y \rcf = \lcf X ,d_\ast{}_{X_0}Y \rcf _{\phi
_0}-\lcf Y ,d_\ast{}_{X_0}X \rcf _{\phi _0},\\[4pt] ({\cal L}_\ast
{}_{X_0})_{\phi _0}P+({\cal L}_{\phi _0})_{X_0}P= 0,
\end{array}
\end{equation}
for $X ,Y \in \Gamma (A)$ and $P\in \Gamma (\wedge ^pA)$, where
$d_\ast{}_{X_0}$ (respectively, ${\cal L}_\ast{}_{X_0}$) is the
$X_0$-differential (respectively, the $X_0$-Lie derivative) of
$A^\ast$. Note that the second equality in (\ref{condcomp}) holds
if and only if
\begin{equation}\label{condcomp2}
\begin{array}{l}
\phi _0(X_0)=0,\qquad \rho (X_0)=-\rho _\ast (\phi _0)\\ ({\cal
L}_\ast {}_{X_0})_{\phi _0}X+\lcf X_0,X\rcf = 0,
\end{array}
\end{equation}
for $X\in \Gamma (A)$ (see \cite{IM}). Very recently, an
interesting characterization of generalized Lie bialgebroids has
been obtained by Grabowski and Marmo \cite{GM} as follows. If we
consider the bracket $\lcf \, ,\, \rcf '_{\phi _0}$ of a
$p$-section $P$ and a $p'$-section $P'$ as the $(p+p'-1)$-section
given by $\lcf P ,P' \rcf '_{\phi _0}=(-1)^{p+1}\lcf P,P'\rcf
_{\phi _0}$ then $((A,\phi _0 ),(A^\ast ,X_0))$ is a generalized
Lie bialgebroid if and only if $d_{\ast X_0}$ is a derivation of
$(\oplus_k\Gamma (\wedge ^kA),\lcf \, ,\, \rcf '_{\phi _0})$, that
is, $$ d_{\ast X_0}\lcf P ,P' \rcf '_{\phi _0}= \lcf d_{\ast X_0}P
,P' \rcf '_{\phi _0}+(-1)^{p+1}\lcf P,d_{\ast X_0}P' \rcf ' _{\phi
_0}$$ for $P\in \Gamma (\wedge ^pA)$ and $P'\in \Gamma (\wedge
^\ast A)$. In the particular case when $\phi _0=0$ and $X_0=0$,
(\ref{condcomp}) is equivalent to the condition $ d_\ast{}\lcf X
,Y \rcf = \lcf X ,d_\ast Y \rcf -\lcf Y ,d_\ast X \rcf .$ Thus,
the pair $((A,0),(A^\ast ,0))$ is a generalized Lie bialgebroid if
and only if the pair $(A,A^\ast)$ is a Lie bialgebroid (see
\cite{K-S,MX}).

On the other hand, if $(M,\Lambda ,E)$ is a Jacobi manifold, then
we proved in \cite{IM} that the pair $\Big ((TM\times \R ,\phi
_0)$,$(T^\ast M \times \R ,X_0)\Big )$ is a generalized Lie
bialgebroid, where $\phi _0$ and $X_0$ are the 1-cocycles on
$TM\times \R$ and $T^\ast M\times \R$ given by $$ \phi _0=(0,1)\in
\Omega ^1(M)\times C^\infty (M,\R )\cong \Gamma (T^\ast M\times
\R), $$ $$ X_0=(-E,0)\in \mathfrak X (M)\times C^\infty (M,\R
)\cong \Gamma (TM\times \R). $$ As a kind of converse, we have the
following result.
\begin{theorem}\cite{IM}\label{Jacobi-base}
Let $((A,\phi _0),\kern-1pt(A^\ast,X_0))$ be a generalized Lie
bialgebroid over $M$. Then, the bracket of functions $\{ \, ,\,
\}_0 :C^\infty(M,\R)\times C^\infty(M,\R)\kern-.75pt\to\kern-.75pt
C^\infty(M,\R)$ given by $$\{ f , g \}_0:=d_{\phi _0}f\cdot
d_\ast{}_{X_0} g, \mbox{ for }f,g\in C^\infty (M,\R ),$$ defines a
Jacobi structure on $M$.
\end{theorem}
If $(\Lambda _0,E_0)$ is the Jacobi structure on $M$ associated
with the Jacobi bracket $\{\, ,\,\} _0$ then
\begin{equation}\label{1.19I}
\# _{\Lambda _0}(\omega _0)=\rho _\ast (\rho ^\ast (\omega
_0)),\quad E_0=\rho _\ast (\phi _0)=- \rho (X_0),
\end{equation}
for $\omega _0\in \Omega ^1(M)$, $\rho ^\ast :\Omega ^1(M)\to
\Gamma (A^\ast )$ being the adjoint operator of the anchor map
$\rho :\Gamma (A)\to \mathfrak X(M)$.

Next, we will recall the construction of the Lie bialgebroid
associated with a generalized Lie bialgebroid (for more details,
see \cite{IM}).

Let $(A,\lcf \, , \, \rcf, \rho )$ be a Lie algebroid over $M$ and
$\phi _0\in \Gamma (A^\ast )$ be a 1-cocycle. Then, there exists
two Lie algebroid structures on the vector bundle $\tilde{A}=
A\times \R\to M\times \R$. First, we consider the map $\ast
:\Gamma (A)\to \mathfrak X (M\times \R)$ given by
\begin{equation}\label{acc-induc}
X^\ast =\rho(X)\circ \pi _1+(\phi _0(X)\circ \pi _1)\frac{\partial
}{\partial t},
\end{equation}
where $\pi _1:M\times\R\to M$ is the canonical projection onto the
first factor. It is easy to prove that $\ast$ is an action of $A$
on $\pi _1$ (see Section \ref{algebroides}). Thus, if $\pi _1^\ast
A$ is the pull-back of $A$ over $\pi _1$ then the vector bundle
$\pi _1^\ast A\to M\times \R$ admits a Lie algebroid structure
$(\lcf \,,\, \rcf \, \bar{ }\, ^{\phi _0}, \bar{\rho}^{\phi _0})$.
It is clear that the vector  bundles $\pi _1^\ast A\to M\times \R$
and $\tilde{A}= A\times \R\to M\times \R$ are isomorphic and that
the space of sections $\Gamma (\tilde{A})$ of $\tilde{A}\to
M\times \R$ can be identified with the set of time-dependent
sections of $A\to M$. Under this identification, the Lie algebroid
structure $(\lcf \, ,\, \rcf \, \bar{ }\, ^{\phi
_0},\bar{\rho}^{\phi _0})$ is given by
\begin{equation}\label{corchbarra}
\lcf \tilde{X},\tilde{Y}\rcf \, \bar{ }\, ^{\phi _0}= \lcf
\tilde{X}, \tilde{Y}\rcf\,\tilde{ } + \phi
_0(\tilde{X})\frac{\partial \tilde{Y}}{\partial t}-\phi
_0(\tilde{Y})\frac{\partial \tilde{X}}{\partial t},\qquad
\bar{\rho}^{\phi _0}(\tilde{X})= \tilde{\rho}(\tilde{X})+\phi
_0(\tilde{X}) \frac{\partial}{\partial t},
\end{equation}
for $\tilde{X},\tilde{Y}\in \Gamma (\tilde{A})$, where $(\lcf \,
,\, \rcf \, \tilde{ },\tilde{\rho})$ is the Lie algebroid
structure on $\pi _1^\ast A$ defined by the zero 1-cocycle and
$\displaystyle \frac{\partial\tilde{X}}{\partial t}$
(respectively, $\displaystyle \frac{\partial\tilde{Y}}{\partial
t}$) denotes the derivative of $\tilde{X}$ (respectively,
$\tilde{Y}$) with respect to the time.

Now, let $\Psi :\tilde{A}\to \tilde{A}$ be the isomorphism of
vector bundles over the identity defined by $\Psi (v,t)=(e^tv,t)$,
for $(v,t)\in A\times \R=\tilde{A}$. Using $\Psi$ and the Lie
algebroid structure $(\lcf \, ,\, \rcf \, \bar{ }\, ^{\phi _0},
\bar{\rho}^{\phi _0})$, one can introduce a new Lie algebroid
structure $(\lcf \, ,\, \rcf \, \hat{ }\, ^{\phi
_0},\hat{\rho}^{\phi _0})$ on the vector bundle $\tilde{A}\to
M\times \R$ in such a way that the Lie algebroids $(\tilde{A},
\lcf \,,\, \rcf \, \bar{ }\, ^{\phi _0}, \bar{\rho}^{\phi _0})$
and $(\tilde{A},\lcf \, ,\, \rcf \, \hat{ }\, ^{\phi
_0},\hat{\rho}^{\phi _0})$ are isomorphic. We have that
\begin{equation}\label{corchtilde}
\begin{array}{c}
\lcf \tilde{X},\tilde{Y}\rcf \, \hat{ }\, ^{\phi _0}=e^{-t}\Big(
\lcf \tilde{X}, \tilde{Y}\rcf\,\tilde{ } + \phi
_0(\tilde{X})(\frac{\partial \tilde{Y}}{\partial
t}-\tilde{Y})-\phi _0(\tilde{Y})(\frac{\partial
\tilde{X}}{\partial t}-\tilde{X})\Big ) ,\\
\\
\hat{\rho}^{\phi _0}(\tilde{X})= e^{-t}\Big( \tilde{\rho}
(\tilde{X})+\phi_0(\tilde{X}) \frac{\partial}{\partial t}\Big ) ,
\end{array}
\end{equation}
for all $\tilde{X},\tilde{Y}\in \Gamma(\tilde{A})$. Moreover, one
may prove the following result.
\begin{theorem}\cite{IM}\label{bialgebrizacion}
Let $((A,\phi _0),(A^\ast ,X_0))$ be a generalized Lie bialgebroid
and $(\Lambda ,E)$ be the induced Jacobi structure on $M$.
Consider on $\tilde{A}=A\times \R$ (resp. $\tilde{A}^\ast=A^\ast
\times \R$) the Lie algebroid structure  $(\lcf \, ,\, \rcf \,
\bar{ }\, ^{\phi _0},\bar{\rho}^{\phi _0})$ (resp. $(\lcf \, ,\,
\rcf _\ast \kern-3pt \hat{ }^{X_0} ,\widehat{\rho_ \ast}^{X_0})$).
Then:
\begin{itemize}
\item[{\it i)}] The pair $(\tilde{A},\tilde{A}^\ast)$ is a Lie
bialgebroid over $M\times \R$.
\item[{\it ii)}] If $\tilde{\Lambda}$ is the induced Poisson
structure on $M\times \R$ then $\tilde{\Lambda}$ is the
Poissonization of the Jacobi structure $(\Lambda ,E)$.
\end{itemize}
\end{theorem}
\setcounter{equation}{0}
\section{Contact groupoids and 1-jet Lie groupoids}\label{caract}
First, we will recall the notion of a contact groupoid.
\begin{definition}\cite{KS}
Let $\gr$ be a Lie groupoid, $\eta \in \Omega ^1(G)$ be a contact
1-form on $G$ and $\sigma :G\to \R$ be an arbitrary function. If
$\oplus _{TG}$ is the partial multiplication in the Lie groupoid
$TG\gpd TM$, we will say that $(\gr ,\eta ,\sigma )$ is a contact
groupoid if and only if
\begin{equation}\label{contgroup}
\eta _{(gh)}(X_g \oplus _{TG} Y_h )=\eta _g(X_g)+\exps \eta
_h(Y_h),\mbox{ for }(X_g,Y_h)\in TG^{(2)}.
\end{equation}
\end{definition}
\begin{remark}
{\rm Actually, the definition of a contact groupoid given in
\cite{KS} is slightly different to the one given here. The
relation between both approaches is the following one. If $(\gr,
\theta, \kappa )$ is a contact groupoid in the sense of \cite{KS}
then $(\gr, \eta ,\sigma )$ is a contact groupoid in the sense of
Definition \ref{contgroup}, where $\sigma (g)=\kappa (g^{-1})$ for
$g\in G$, and $\eta _g$ is the inverse of $\theta _{g^{-1}}$ in
the Lie groupoid $T^\ast G\gpd A^\ast G$. }
\end{remark}
If $(\gr ,\eta ,\sigma )$ is a contact groupoid then, using the
associativity of $\oplus _{TG}$, we deduce that $\mult$ is a
multiplicative function, that is,
\begin{equation}\label{multiplicidad}
\sigma (gh)=\sigma (g)+\sigma (h)
\end{equation}
for $(g,h)\in G^{(2)}$. In particular, $\sigma _{|\epsilon
(M)}\equiv 0$ and therefore, using (\ref{contgroup}), it follows
that $\eta _{\tilde{x}}(\epsilon ^x_\ast (X_x))$ $=0$, for $x\in
M$ and $X_x\in T_xM$. Thus, if $\iota:G\to G$ is the inversion of
$G$, we obtain that $\iota ^\ast \eta =-e^{-\sigma }\eta$. This
implies that $G$ is a contact groupoid in the sense of \cite{D2}.
Using this fact, we deduce the following result.
\begin{proposition}\label{utiles}
Let $(\gr ,\eta ,\sigma )$ be a contact groupoid and suppose that
$dim\, G=2n+1$. Then:
\begin{itemize}
\item[{\it i)}] If $g$ and $h$ are composable elements of $G$, we
have that
\begin{equation}\label{2.2'}
\begin{array}{lll} (\delta \eta )_{gh} (X_g\oplus
_{TG}Y_h,X_g'\oplus _{TG}Y'_h)&=& (\delta \eta) _g
(X_g,X_g')+e^{\sigma (g)}(\delta \eta )_h(Y_h,Y'_h)
\\&&+e^{\sigma (g)}(X_g(\sigma )\eta _h(Y'_h)- X'_g(\sigma )\eta
_h(Y_h)),\end{array}
\end{equation}
for $(X_g,Y_h),(X'_g,Y'_h)\in TG^{(2)}$.
\item[{\it ii)}] $M\cong \epsilon (M)$ is a Legendre submanifold of $G$,
that is, $\epsilon ^\ast \eta =0$ and $dim\,\epsilon (M)=dim\,
M=n$.
\item[{\it iii)}] If $(\Jacobi )$ is the Jacobi structure associated with
the contact 1-form $\eta$, then $E$ is a right-invariant vector
field on $G$ and $E(\sigma )=0$. Moreover, if $X_0\in \Gamma (AG)$
is the section of the Lie algebroid $AG$ of $G$ satisfying
$E=-\der{X_0}$, we have that
\begin{equation}\label{2.2''}
\sostP (\delta \sigma )=\der{X_0}-e^{-\sigma}\izq{X_0}.
\end{equation}
\vspace{-20pt}

\item[{\it iv)}] If $\s ^T$, $\t ^T$ and $\epsilon ^T$ (respectively,
$\tilde{\s}$, $\tilde{\t}$ and $\tilde{\epsilon}$) are the
projections and the inclusion in the Lie groupoid $TG\gpd TM$
(respectively, $T^\ast G \gpd A^\ast G$) then,
$$e^{-\sigma}\sostP\circ \tilde{\epsilon} \circ
\tilde{\s}=\epsilon ^T\circ \s ^T\circ \sostP,\qquad \sostP\circ
\tilde{\epsilon} \circ \tilde{\t}= \epsilon ^T\circ \t ^T\circ
\sostP .$$
\end{itemize}
\end{proposition}
\prueba Using the results in \cite{D2}, we directly deduce {\it
i)}, {\it ii)} and {\it iii)}.

Now, we will prove {\it iv)}. Suppose that $\omega _g\in T^\ast
_gG$. Then, from {\it ii)}, we conclude that $$\eta
_{\widetilde{\s (g)}}(\expsn \sostP
(\tilde{\epsilon}(\tilde{\s}(\omega _g)))=\eta _{\widetilde{\s
(g)}}(\epsilon _\ast ^{\s (g)}(\s _\ast ^g(\sostP (\omega
_g)))=0.$$ Furthermore, if $X_{\s (g)}\in A_{\s (g)}G$, it follows
that $$\epsilon _\ast ^{\s (g)}(\s _\ast ^g(\sostP (\omega
_g)))=\iota ^g_\ast (\sostP (\omega _g))\oplus _{TG}\sostP (\omega
_g),\qquad X_{\s(g)}=0_{T_{g\inv} G}\oplus _{TG}(L_g)_\ast
^{\widetilde{\s(g)}}(X_{\s(g)})$$ and consequently, using
(\ref{T*G}), (\ref{contgroup}), (\ref{2.2'}) and (\ref{2.2''}) and
the fact that $\sigma$ is a multiplicative function, we obtain
that $$(\delta \eta )_{\widetilde{\s (g)}}(\epsilon _\ast ^{\s
(g)}(\s _\ast ^g(\sostP (\omega _g))),X_{\s (g)})=(\delta \eta
)_{\widetilde{\s (g)}}(\expsn \sostP
(\tilde{\epsilon}(\tilde{\s}(\omega _g))),X_{\s (g)}).$$ On the
other hand, from (\ref{T*G}) and {\it ii)}, we deduce that
$$(\delta \eta )_{\widetilde{\s (g)}}(\epsilon _\ast ^{\s (g)}(\s
_\ast ^g(\sostP (\omega _g))),\epsilon _\ast ^{\s(g)}(Y_{\s
(g)}))=(\delta \eta )_{\widetilde{\s (g)}}(\expsn \sostP
(\tilde{\epsilon}( \tilde{\s}(\omega _g))),\epsilon _\ast
^{\s(g)}(Y_{\s (g)}))=0$$ for $Y_{\s(g)}\in T_{\s(g)}M$.

The above facts imply that $\epsilon ^T (\s ^T (\sostP (\omega
_g)))=\expsn \sostP (\tilde{\epsilon }(\tilde{\s}(\omega _g)))$.
In a similar way, one may prove that
$\sostP(\tilde{\epsilon}(\tilde{\t}(\omega _g)))=\epsilon ^T(\t
^T(\sostP (\omega _g)))$. \QED

Using again the results in \cite{D2}, we have that
\begin{proposition}\label{utiles2}
Let $(\gr ,\eta ,\sigma )$ be a contact groupoid and $\mathfrak X
_L(G)$ be the set of left-invariant vector fields on $G$. Denote
by $(\Jacobi )$ the Jacobi structure on $G$ associated with the
contact 1-form $\eta$, by $X_0\in \Gamma (AG)$ the section of the
Lie algebroid $AG$ of $G$ satisfying $E=-\der{X_0}$ and by ${\cal
I}:\Omega ^1(M)\times C^\infty (M,\R)\to \mathfrak X (G)$ the map
defined by
\begin{equation}\label{identificacion}
{\cal I}(\omega _0,f_0)=\sostP (e^{\sigma}\s ^\ast \omega _0)-(\s
^\ast f_0)\izq{X_0}.
\end{equation}
Then:
\begin{itemize}
\item[{\it i)}] ${\cal I}$ defines an isomorphism of $C^\infty (M,\R
)$ modules between the spaces $\Omega ^1(M)\times C^\infty (M,$ $
\R)$ and  $\mathfrak X _L(G)$.
\item[{\it ii)}] The base manifold $M$ admits a Jacobi structure
$(\Lambda _0,E_0)$ in such a way that the projection $\t$ is a
Jacobi antimorphism and the pair $(\s ,e^\sigma )$ is a conformal
Jacobi morphism, that is,
\begin{equation}\label{Jacobi-base-contacto}
\begin{array}{l} \Lambda _0(\s
(g))=\exps \s _\ast ^g(\Lambda (g)),\qquad E_0(\s (g))=\s _\ast
^g(X_{e^\sigma}(g)),\\\Lambda _0(\t (g))=- \t _\ast ^g(\Lambda
(g)),\qquad E_0(\t (g))=-\t _\ast ^g(E(g)),
\end{array}
\end{equation}
for all $g\in G$, where $X_{e^\sigma}=e^{\sigma}\sostP (\delta
\sigma )+e^\sigma E$ is the hamiltonian vector field of the
function $e^\sigma$ with respect to the Jacobi structure $(\Jacobi
)$.
\item[{\it iii)}] The map ${\cal I}$ induces an isomorphism between
the Lie algebroids $(T^\ast M\times\R, \lcf \, ,\, \rcf _{(\Lambda
_0,E_0)},$ $\widetilde{\#}_{(\Lambda _0,E_0)})$ and $AG$.
\end{itemize}
\end{proposition}
\begin{remark}{\rm Denote also by ${\cal I}: T^\ast M\times\R\to
AG$ the Lie algebroid isomorphism induced by the isomorphism of
$C^\infty (M,\R)$-modules ${\cal I}:\Omega ^1(M)\times C^\infty
(M,\R)\to \mathfrak X _L(G)$. Then, from (\ref{identificacion})
and since $\sigma$ is a multiplicative function, it follows that
\begin{equation}\label{2.2'''}
{\cal I}(\omega _x,\gamma )=\sostP ((\s _\ast ^{\tilde{x}})^\ast
(\omega _x)) - \gamma X_0(x),
\end{equation}
for $(\omega _x,\gamma )\in T^\ast _xM\times\R$.}
\end{remark}
Now, let $\gr$ be a Lie groupoid and  $\mult$ be a multiplicative
function. Then, there exists a natural right action of the tangent
groupoid $TG\gpd TM$ on the projection $\pi _1:TM\times \R \to TM$
given by $$ (v_x,\lambda )\cdot X_g=(v_x ,X_g(\sigma )+\lambda),$$
for $(v_x,\lambda )\in TM\times\R$ and $X_g\in T_gG$ satisfying
$\t ^T(X_g)=\pi _1(v_x, \lambda )$ (see Examples
\ref{ej-grupoides}). The resulting action groupoid is isomorphic
to $TG\times \R\gpd TM\times \R$ with projections $(\s
^T)_\sigma$, $(\t ^T)_\sigma$, partial multiplication $\oplus
_{TG\times \Rp}$, inclusion $(\epsilon ^T)_\sigma$ and inversion
$(\iota ^T)_\sigma$ given by
\begin{equation}\label{TGR}
\begin{array}{l} (\s ^T)_\sigma (X_g,\lambda )=(\s ^T(X_g),
X_g(\sigma )+\lambda ),\mbox{ for }(X_g,\lambda )\in T_gG\times\R,
\\(\t ^T)_\sigma (Y_h,\mu )=(\t ^T(Y_h),\mu ),
\mbox{ for }(Y_h,\mu )\in T_hG\times\R ,\\ (X_g,\lambda )\oplus
_{TG\times \Rp}(Y_h,\mu )=(X_g\oplus_{TG} Y_h,\lambda ), \mbox{ if
}(\s ^T)_\sigma (X_g,\lambda )=(\t ^T)_\sigma (Y_h,\mu ),\\
(\epsilon ^T)_\sigma (X_x,\lambda )=(\epsilon ^T(X_x),\lambda
),\mbox{ for }(X_x,\lambda )\in T_xM\times\R ,\\ (\iota ^T)_\sigma
(X_g,\lambda )=(\iota ^T(X_g),X_g(\sigma )+\lambda ),\mbox{ for
}(X_g,\lambda )\in T_gG\times\R .
\end{array}
\end{equation}
Now, suppose that $(\gr ,\eta ,\sigma )$ is a contact groupoid.
Since $\eta$ is a contact 1-form, the map $\# _{(\delta \eta ,\eta
)}:T G\times \R\to T^\ast G\times \R$ given by
\begin{equation}\label{morf-contacto}
\# _{(\delta \eta ,\eta )}(X_g,\lambda )=(-i(X_g)(\delta \eta)_g
-\lambda \,\eta _g,\eta _g(X_g))
\end{equation}
is an isomorphism of vector bundles. The inverse map of $\#
_{(\delta \eta ,\eta )}$ is the homomorphism $\sostJ :T^\ast
G\times \R\to TG\times \R$ defined by
\begin{equation}\label{morf-Jacobi}
\sostJ (\omega_g,\gamma )=(\sostP(\omega _g)+\gamma \,E(g),-\omega
_g(E(g))),
\end{equation}
where $(\Jacobi )$ is the Jacobi structure associated with the
contact 1-form $\eta$.

On the other hand, if $A^\ast G$ is the dual bundle to the Lie
algebroid $AG$ then, since $\epsilon (M)$ is a Legendre
submanifold of $G$, the map $\psi _0: TM\times\R\to A^\ast G$
given by
\begin{equation}\label{2.3'''}
\psi _0(X_x,\lambda )=(-i(\epsilon _\ast ^x(X_x))(\delta \eta
)_{\tilde{x}}-\lambda\,\eta _{\tilde{x}})_{|A_xG},\mbox{ for
}(X_x,\lambda )\in T_xM\times \R
\end{equation}
is an isomorphism of vector bundles. Note that $\# _{(\delta \eta
,\eta )}(\epsilon _\ast ^x(X_x),\lambda )=(\tilde{\epsilon}(\psi
_0 (X_x,\lambda )),0)$ and thus the inverse map $\varphi _0:A^\ast
G\to TM\times\R$ of $\psi _0$ is defined by
\begin{equation}\label{2.3iv}
\varphi _0(\omega _x)=(\s _\ast ^{\tilde{x}}(\sostP
(\tilde{\epsilon}(\omega _x))),-\omega _x(E_{\tilde{x}}-\epsilon
_\ast ^x(\t _\ast ^{\tilde{x}}(E_{\tilde{x}})))),
\end{equation}
$\tilde{\epsilon}:A^\ast G\to T^\ast G$ being the inclusion of
identities in the Lie groupoid $T^\ast G\gpd A^\ast G$.

Next, we consider the maps $\tilde{\s}_\sigma, \tilde{\t}_\sigma
:T^\ast G\times\R \to A^\ast G$, $\tilde{\epsilon}_\sigma :A^\ast
G\to T^\ast G\times\R$ and $\tilde{\iota}_\sigma :T^\ast
G\times\R\to T^\ast G\times \R$ given by
\begin{equation}\label{isomorfismos}
\begin{array}{l} \tilde{\s}_\sigma =\psi _0\circ (\s ^T
)_\sigma \circ \sostJ ,\quad \tilde{\t}_\sigma =\psi _0\circ (\t
^T )_\sigma \circ \sostJ,\\ \tilde{\epsilon}_\sigma =\# _{(\delta
\eta ,\eta )}\circ (\epsilon ^T )_\sigma \circ \varphi _0,\quad
\tilde{\iota}_\sigma =\# _{(\delta \eta ,\eta )}\circ (\iota
^T)_\sigma \circ \sostJ
\end{array}
\end{equation}
and the partial multiplication $\oplus _{T^\ast G\times\Rp }$
defined as follows. If $(\omega _g,\gamma ),(\nu _h,\zeta )\in
T^\ast G\times \R$ satisfy $\tilde{\s}_\sigma (\omega _g,\gamma
)=\tilde{\t}_\sigma (\nu _h,\zeta )$ then
\begin{equation}\label{morf-grupos}
(\omega _g,\gamma )\oplus _{T^\ast G\times \Rp}(\nu _h,\zeta )=\#
_{(\delta \eta ,\eta )}\Big (\sostJ (\omega _g,\gamma )\oplus _{T
G\times \Rp}\sostJ (\nu _h,\zeta )\Big ) .
\end{equation}
It is clear $\tilde{\s}_\sigma$, $\tilde{\t}_\sigma$,
$\tilde{\epsilon}_\sigma$, $\tilde{\iota}_\sigma$ and the partial
multiplication $\oplus _{T^\ast G\times\Rp }$ are the structural
functions of a Lie groupoid structure  in $T^\ast G\times\R$ over
$A^ \ast G$. In addition, the map $\sostJ :T^\ast G\times\R \to
TG\times \R$ is a Lie groupoid isomorphism over $\varphi _0
:A^\ast G\to TM\times\R$.
\begin{lemma}
If $\tilde{\s}$, $\tilde{\t}$, $\oplus _{T^\ast G}$,
$\tilde{\epsilon}$ and $\tilde{\iota}$ are the structural
functions of the Lie groupoid $T^\ast G\gpd A^\ast G$, we have
that
\begin{equation}\label{T*GR}
\begin{array}{l}
\tilde{\s}_\sigma(\omega _g,\gamma )=e^{-\sigma(g)}
\tilde{\s}(\omega _g), \mbox{ for }(\omega _g,\gamma )\in T_g^\ast
G\times \R,\\ \tilde{\t}_\sigma(\nu _h,\zeta )=\tilde{\t}(\nu
_h)-\zeta \, (\delta \sigma )_{\widetilde{\t (h)}}{}_{|A_{\t
(h)}G} ,\mbox{ for }(\nu _h,\zeta )\in T_h^\ast G\times \R,\\ \Big
( (\omega _g,\gamma )\oplus _{T^\ast G\times \Rp}(\nu _h,\zeta
)\Big )=\Big ( (\omega _g +e^{\sigma (g)}\zeta \,(\delta \sigma )
_g)\oplus _{T^\ast G}(e^{\sigma (g)}\nu _h),\gamma +e^{\sigma
(g)}\zeta \Big )\\ \tilde{\epsilon}_\sigma (\omega
_x)=(\tilde{\epsilon}(\omega _x),0),\mbox{ for }\omega _x\in
A^\ast _xG,\\ \tilde{\iota}_\sigma (\omega _g,\gamma )=(e^{-\sigma
(g)}(\tilde{\iota}(\omega _g)-\gamma (\delta \sigma
)_{g\inv}),-e^{-\sigma (g)}\gamma ),\mbox{ for }(\omega _g,\gamma
)\in T^\ast _gG\times\R .
\end{array}
\end{equation}
\end{lemma}
\prueba A long computation, using (\ref{T*G}), (\ref{contgroup}),
(\ref{multiplicidad}),  (\ref{TGR})-(\ref{morf-grupos}) and
Proposition \ref{utiles}, proves the result.\QED

Note that the maps $\tilde{\s}_\sigma$, $\tilde{\t}_\sigma$,
$\tilde{\epsilon}_\sigma$, $\tilde{\iota}_\sigma$ and the partial
multiplication $\oplus _{T^\ast G\times\Rp }$ do not depend on the
contact 1-form $\eta$. In fact, one may prove the following
result.
\begin{theorem}\label{contacto}
Let $\gr$ be an arbitrary Lie groupoid with Lie algebroid $AG$ and
$\sigma :G\to\R$ be a multiplicative function. Then:
\begin{itemize}
\item[{\it i)}] The product manifold $T^\ast G\times\R$ admits a
Lie groupoid structure over $A^\ast G$ with structural functions
given by (\ref{T*GR}).
\item[{\it ii)}] If $\eta _G$ is the canonical contact 1-form on
$T^\ast G\times\R$ and $\bar{\pi}_G:T^\ast G\times \R\to G$ is the
canonical projection then $(T^\ast G\times \R\gpd A^\ast G,\eta
_G, \sigma \circ \bar{\pi }_G)$ is a contact groupoid.
\end{itemize}
\end{theorem}
\prueba Since $\sigma$ is a multiplicative function, we obtain
that
\begin{equation}\label{2.13I}
\epsilon ^\ast \sigma =0.
\end{equation}
Moreover, if $(g,h)\in G^{(2)}$ and $\s(g)=\t (h)=x\in M$ then,
from (\ref{T*G}), it follows that
\begin{equation}\label{2.13II}
\begin{array}{l}
\tilde{\s}((\delta \sigma )_g)=\tilde{\t}((\delta \sigma
)_h)=(\delta \sigma )_{\tilde{x}}{}_{|A_xG},\\ (\delta \sigma
)_{gh}=(\delta \sigma )_g\oplus_{T^\ast G} (\delta \sigma )_h .
\end{array}
\end{equation}
In addition, using again (\ref{T*G}) and the fact that $\sigma$ is
a multiplicative function, we have that
\begin{equation}\label{2.13III}
\tilde{\epsilon}((\delta \sigma )_{\tilde{x}}{}_{|A_xG})=(\delta
\sigma )_{\tilde{x}},\quad \tilde{\iota} ((\delta \sigma
)_g)=(\delta \sigma )_{g\inv},
\end{equation}
for $x\in M$ and $g\in G$.

Thus, from (\ref{T*GR})-(\ref{2.13III}), we deduce {\it i)}.

Now, let  $G\times\R \gpd M$ be the semi-direct Lie groupoid with
projections $\s '$, $\t '$, partial multiplication $m'$, inclusion
$\epsilon '$ and inversion $\iota '$ defined by
\begin{equation}\label{2.14}
\begin{array}{l}
\s '(g,\gamma )=\s (g),\quad \t '(g,\gamma )=\t (g), \mbox{ for
}(g,\gamma )\in G\times \R ,\\ m'((g,\gamma ),(h,\zeta ))=
(gh,\gamma +\exps \zeta) , \mbox{ for } ((g,\gamma ),(h,\zeta
))\in (G\times \R )^{(2)},\\ \epsilon '(x)=(\epsilon (x),0),\mbox{
for }x\in M,\\ \iota ' (g,\gamma )=(\iota (g),-\expsn \gamma ),
\mbox{ for }(g,\gamma )\in G\times\R.\end{array}
\end{equation}
Using (\ref{2.14}), one may prove that the partial multiplication
$\oplus _{T(G\times \Rp )}$ in the tangent Lie groupoid
$T(G\times\R)\gpd TM$ is given by
\begin{equation}\label{e1}
\begin{array}{rcl}
\Big ( X_g+\psi \frac{\partial}{\partial t}_{|\gamma }\Big )\oplus
_{T(G\times \Rp )} \Big ( Y_h+\varphi \frac{\partial}{\partial
t}_{|\zeta }\Big ) &=&(X_g\oplus _{TG}Y_h)\\&&\kern-15pt +(\psi
+\exps (\zeta X_g (\sigma )+\varphi ))\frac{\partial}{\partial
t}_{|\gamma +\exps \zeta }.\end{array}
\end{equation}
Next, we consider the map $\tilde{\pi}_G:T^\ast G\times \R\to
G\times\R$ given by $$\tilde{\pi}_G (\omega _g,\gamma )=(\pi
_G(\omega _g),\gamma ),\mbox{ for }(\omega _g,\gamma )\in T^\ast
_gG\times \R,$$ where $\pi _G:T^\ast G\to G$ is the canonical
projection. From (\ref{T*GR}) and (\ref{2.14}), we deduce that
$\tilde{\pi}_G$ is a Lie groupoid morphism over the map
$\tilde{\pi}_0: A^\ast G\to M$ defined by $$\tilde{\pi }_0 (\omega
_x)=x,\mbox{ for }\omega _x\in A^\ast _xG.$$ Therefore, the
tangent map to $\tilde{\pi}_G$, $T\tilde{\pi}_G :T(T^\ast
G\times\R)\to T(G\times\R )$, given by
\begin{equation}\label{2.15I}
T\tilde{\pi}_G( X_{\omega _g}+\psi \frac{\partial}{\partial
t}_{|\gamma })= (\pi _G)_\ast ^{\omega _g}(X_{\omega _g})+\psi
\frac{\partial}{\partial t}_{|\gamma },
\end{equation}
for $X_{\omega _g}+\psi \frac{\partial}{\partial t}_{|\gamma } \in
T_{(\omega _g,\gamma )}(G\times \R)$, is also a Lie groupoid
morphism (over the map $T\tilde{\pi }_0:T(A^\ast G)\to TM$)
between the tangent Lie groupoids $T(T^\ast G\times\R )\gpd
T(A^\ast G)$ and $T(G\times\R)\gpd TM$.

On the other hand, if $\eta _G$ is the canonical contact 1-form on
$T^\ast G\times \R$ then $\eta _G=\lambda _G-\delta t$, $\lambda
_G$ being the Liouville 1-form on $T^\ast G$, and (see
(\ref{2.15I}))
\begin{equation}\label{e2}
\begin{array}{lll} \eta _G{}_{(\omega _g,\lambda
)}(X_{\omega _g}+\psi \frac{\partial}{\partial t}_{|\gamma
})&=&\lambda _G{}_{(\omega _g )}(X_{\omega _g})-\delta t_{|\gamma
}(\psi \frac{\partial}{\partial t}_{|\gamma })=\omega _g((\pi
_G)_\ast ^{\omega _g}(X_{\omega _g}))-\psi \\&=&(\omega _g-\delta
t_{|\gamma })(T\tilde{\pi}_G (X_{\omega _g}+\psi
\frac{\partial}{\partial t}_{|\gamma })).\end{array}
\end{equation}
Thus, using (\ref{T*GR}), (\ref{e1}), (\ref{2.15I}), (\ref{e2})
and the fact that $T\tilde{\pi}_G$ is a Lie groupoid morphism, we
conclude that $$\eta _G{}_{((\omega _g,\gamma )\oplus _{T^\ast
G\times \Rp}(\nu _h,\zeta ))}=\eta _G{}_{(\omega _g,\gamma
)}\oplus _{T^\ast (T^\ast G\times \Rp)} ( \exps \eta _G{}_{(\nu
_h,\zeta )}),$$ that is, $(T^\ast G\times \R\gpd A^\ast G,\eta
_G,\bar{\sigma })$ is a contact groupoid, where $\bar{\sigma}\in
C^\infty (T^\ast G\times\R )$ is the function given by
$\bar{\sigma}=\sigma\circ \bar{\pi}_G$.\QED
\begin{remark}\label{recup}
{\rm Let $\gr$ be a Lie groupoid, $\mult$ be a multiplicative
function and $TG\times\R\gpd TM\times\R$, $T^\ast G\times\R\gpd
A^\ast G$ be the corresponding Lie groupoids with structural
functions given by (\ref{TGR}) and (\ref{T*GR}). If $\sigma$
identically vanishes then we recover, by projection, the Lie
groupoids $TG\gpd TM$ and $T^\ast G\gpd A^\ast G$ (see Examples
\ref{ej-grupoides}).}
\end{remark}
\begin{remark}\label{2.8'}
{\rm {\it i)} A Lie groupoid $\gr$ is said to be symplectic if $G$
admits a symplectic 2-form $\Omega$ in such a way that the graph
of the partial multiplication in $G$ is a Lagrangian submanifold
of the symplectic manifold $(G\times G\times G,\Omega \oplus
\Omega \oplus (-\Omega ))$ (see \cite{CDW}). If $\gr$ is an
arbitrary Lie groupoid with Lie algebroid $AG$ and on the
cotangent Lie groupoid $T^\ast G$ we consider the canonical
symplectic $2$-form  $\Omega _G=-\delta \lambda _G$ then $T^\ast
G$ is a symplectic groupoid over $A^\ast G$ (see \cite{CDW}).

{\it ii)} Let $\gr$ be a symplectic groupoid with exact symplectic
2-form $\Omega =-\delta \theta$. Then, since $\R$ is a Lie group,
the product manifold $G\times\R$ is a Lie groupoid over $M$ (see
Examples \ref{ej-grupoides}, 3). In addition, $(G\times\R\gpd
M,\eta ,0)$ is a contact groupoid, where $\eta$ is the 1-form on
$G\times\R$ given by $\eta =\pi _1^\ast (\theta )-\pi _2^\ast
(\delta t)$, and $\pi _1:G\times\R\to G$, $\pi _2:G\times\R\to \R$
are the canonical projections (see \cite{Lib2}). In particular, if
$\gr$ is an arbitrary Lie groupoid with Lie algebroid $AG$ then we
have that $(T^\ast G\times\R\gpd A^\ast G, \eta _G,0)$ is a
contact groupoid, $\eta _G$ being the canonical contact 1-form on
$T^\ast G\times\R$. Note that, using Theorem \ref{contacto}, we
directly deduce this result. }
\end{remark}
Let $\gr$ be an arbitrary Lie groupoid with Lie algebroid $AG$ and
$\mult$ be a multiplicative function. From Proposition
\ref{utiles2}, it follows that the contact groupoid structure on
$T^\ast G\times \R$ induces a Jacobi structure on the vector
bundle $A^\ast G$. Next, we will describe such a Jacobi structure.
For this purpose, we will recall the definition of the linear
Jacobi structure associated with a Lie algebroid and a 1-cocycle
on it (for more details, see \cite{IM0}).

Suppose that $(L,\lcf \, ,\, \rcf ,\rho )$ is a Lie algebroid over
$M$ and denote by $\Lambda_{L^\ast }$ the corresponding linear
Poisson structure on $L^\ast$ (see Section \ref{algebroides}). If
$\omega _0\in \Gamma (L^\ast)$ is a 1-cocycle of $L$, $\Delta$ is
the Liouville vector field of $L^\ast$ and $\omega _0^v\in
\mathfrak X (L^\ast)$ is the vertical lift of $\omega _0$, we have
that the pair $(\Lambda_{(L^\ast ,\omega _0)},E_{(L^\ast ,\omega
_0)})$ is a Jacobi structure on $L^\ast$, where
\begin{equation}\label{Jacobi-lineal}
\Lambda_{(L^\ast ,\omega _0)}=\Lambda_{L^\ast }+\Delta \wedge
\omega _0^v, \qquad E_{(L^\ast ,\omega _0)}=-\omega _0^v.
\end{equation}
The Jacobi bracket $\{\, ,\,\} _{(L^\ast ,\omega _0)}$ associated
with the Jacobi structure $(\Lambda_{(L^\ast ,\omega
_0)},E_{(L^\ast ,\omega _0)})$ is characterized by the following
conditions
\begin{equation}\label{2.16II}
\{ \tilde{X},\tilde{Y} \}_{(L^\ast ,\omega _0)} =\widetilde{\lcf
X,Y\rcf },\quad \{ \tilde{X}, 1 \}_{(L^\ast ,\omega _0)} =\omega
_0(X)\circ \tau ^\ast ,
\end{equation}
for $X,Y\in \Gamma (L)$, $\tau ^\ast :L^\ast \to M$ being the
bundle projection. Here, if $Z\in \Gamma (L)$, we denote by
$\tilde{Z}$ the corresponding linear function on $L^\ast$ (see
\cite{IM0}).
\begin{theorem}\label{2.8}
Let $\gr$ be a Lie groupoid with Lie algebroid $AG$ and $\mult$ be
a multiplicative function. If $\bar{\pi}_G:T^\ast G\times\R\to G$
is the canonical projection, $\eta _G$ is the canonical contact
1-form on $T^\ast G\times\R$ and $(\Lambda _0,E_0)$ is the Jacobi
structure on $A^\ast G$ induced by the contact groupoid $(T^\ast
G\times\R\gpd A^\ast G,\eta _G,\bar{\sigma}=\sigma \circ
\bar{\pi}_G)$ then
\begin{equation}\label{base=lineal}
\Lambda _0=\Lambda_{(A^\ast G,\phi _0)},\quad E_0=E_{(A^\ast G
,\phi_0)},
\end{equation}
where $\phi _0\in \Gamma (A^\ast G)$ is the 1-cocycle of the Lie
algebroid $AG$ defined by
\begin{equation}\label{2.16III}
\phi _0(x)(X_x)=X_x(\sigma ), \mbox{ for }x\in M\mbox{ and }X_x\in
A_xG.
\end{equation}
\end{theorem}
\prueba Denote by $\pi _1: T^\ast G\times\R\to T^\ast G$ the
canonical projection onto the first factor. It is easy to prove
that $\pi _1$ is a Jacobi morphism between the contact manifold
$(T^\ast G\times\R, \eta _G)$ and the symplectic manifold $(T^\ast
G,\Omega _G)$. This means that
\begin{equation}\label{2.17}
\{ f\circ \pi _1,g \circ \pi _1 \}_{T^\ast G\times\Rp}=\{ f,g
\}_{T^\ast G}\circ \pi _1,
\end{equation}
for $f,g\in C^\infty (T^\ast G,\R )$,  $\{ \, ,\, \}_{T^\ast
G\times\Rp}$ (respectively, $\{ \, ,\, \}_{T^\ast G}$) being the
Jacobi bracket (respectively, Poisson bracket) associated with the
contact 1-form $\eta _G$ (respectively, the symplectic 2-form
$\Omega _G$).

Now, suppose that $\{\, ,\, \}_0$ is the Jacobi bracket associated
with the Jacobi structure $(\Lambda _0,E_0)$. From
(\ref{Jacobi-base-contacto}), it follows that
\begin{equation}\label{2.17I}
\tilde{\s}_\sigma ^\ast \{ \tilde{f},\tilde{g} \} _0 =
e^{-\bar{\sigma}} \{ e^{\bar{\sigma}} \tilde{\s}_\sigma ^\ast
\tilde{f} ,e^{\bar{\sigma}}\tilde{\s}_\sigma ^\ast \tilde{g}
\}_{T^\ast G\times\Rp}
\end{equation}
for $\tilde{f},\tilde{g}\in C^\infty (T^\ast G\times\R,\R )$.
Thus, if $X,Y\in \Gamma (AG)$ and $\tilde{X},\tilde{Y}$ are the
corresponding linear functions on $A^\ast G$, then (see
(\ref{T*GR}), (\ref{2.17}) and (\ref{2.17I}))
\begin{equation}\label{2.17II}
\begin{array}{lll}
\{ \tilde{X},\tilde{Y} \} _0(\tilde{\s}_\sigma (\omega _g,\gamma
))&=& (e^{-\bar{\sigma}} \{ \tilde{\s} ^\ast (\tilde{X}) \circ \pi
_1 ,\tilde{\s}^\ast (\tilde{Y})\circ \pi _1 \}_{T^\ast
G\times\Rp})(\omega _g,\gamma )\\ &=&
e^{-\sigma(g)}\{\tilde{\s}^*(\tilde{X}),
\tilde{\s}^*(\tilde{Y})\}_{T^*G}(\omega_g),
\end{array}
\end{equation}
for $(\omega _g,\gamma )\in T^\ast _gG\times \R$. On the other
hand, using the results in \cite{CDW}, we have that
\begin{equation}\label{proyeccion}
(\pi _G)^{\nu _h}_\ast (X^{\Omega _G}_{\tilde{\s} ^\ast
(\tilde{X})}(\nu _h ))=\izq{X}(h),\quad (\pi _G)^{\nu _h}_\ast
(X^{\Omega _G}_{\tilde{\s} ^\ast (\tilde{Y})}(\nu _h
))=\izq{Y}(h),
\end{equation}
for $h\in G$ and $\nu _h \in T^\ast _hG$, where $X^{\Omega
_G}_{\tilde{\s} ^\ast (\tilde{X})}$ (respectively, $X^{\Omega
_G}_{\tilde{\s} ^\ast (\tilde{Y})}$) is the hamiltonian vector
field of the function $\tilde{\s} ^\ast (\tilde{X})$
(respectively, $\tilde{\s} ^\ast (\tilde{Y})$) with respect to the
symplectic structure $\Omega _G$.
 Therefore, ${\cal L}_{X^{\Omega
_G}_{\tilde{\s} ^\ast (\tilde{X})}}\lambda_G={\cal L}_{X^{\Omega
_G}_{\tilde{\s} ^\ast (\tilde{Y})}}\lambda_G=0$ and from
(\ref{2.17II}) and (\ref{proyeccion}), we conclude that
$$\begin{array}{rcl} \{ \tilde{X},\tilde{Y} \} _0
(\tilde{\s}_\sigma (\omega_g,\gamma ))&=&\expsn \lambda
_G(\omega_g)[X^{\Omega _G}_{\tilde{\s} ^\ast (\tilde{X})},
X^{\Omega _G}_{\tilde{\s} ^\ast (\tilde{Y})}](\omega_g)\\&=&
\expsn \omega _g(\izq{\lcf X, Y\rcf }(g))=\tilde{\s}_\sigma
(\omega _g,\gamma )(\lcf X,Y\rcf (\s (g))),\end{array}$$
$\lcf\;,\;\rcf$ being the Lie bracket on $AG$. Consequently,
\begin{equation}\label{lineal1}
\{ \tilde{X},\tilde{Y} \} _0=\widetilde{\lcf X,Y\rcf }.
\end{equation}
Next, we will show that
\begin{equation}\label{lineal2}
\{\tilde{X},1\}_0=\phi_0(X)\circ \tau_{A^*G},
\end{equation}
where $\tau_{A^\ast G}:A^\ast G\to M$ is the bundle projection.
Using (\ref{T*GR}), (\ref{2.17}) and (\ref{2.17I}), it follows
that $$\begin{array}{rcl} \{ \tilde{X},1 \} _0 (\tilde{\s}_\sigma
(\omega _g,\gamma ))&=& (e^{-\bar{\sigma}} \{ \tilde{\s} ^\ast
(\tilde{X})\circ\pi_1 ,e^{{\sigma}\circ \pi_G}\circ \pi_1
\}_{T^\ast G\times\Rp})(\omega _g,\gamma )\\&=& e^{-{\sigma}(g)}
\{ \tilde{\s} ^\ast (\tilde{X}) ,e^{{\sigma}\circ \pi_G}
\}_{T^\ast G}(\omega _g )= (\pi _G)^{\omega _g}_\ast (X^{\Omega
_G}_{\tilde{\s} ^\ast (\tilde{X})}(\omega _g))(\sigma ).
\end{array}$$

Thus, from (\ref{2.16III}) and (\ref{proyeccion}), we obtain that
\[
\{ \tilde{X},1 \} _0 (\tilde{\s}_\sigma (\omega _g,\gamma ))=
(\phi _0(X)\circ \tau _{A^\ast G})(\tilde{\s}_\sigma (\omega
_g,\gamma )).\] This implies that (\ref{lineal2}) holds.

Finally, using (\ref{lineal1}) and (\ref{lineal2}), we deduce
(\ref{base=lineal}). \QED

\section{Jacobi groupoids}
\setcounter{equation}{0}
\subsection{Jacobi groupoids: definition and examples} Motivated
by the results obtained in Section \ref{caract} about contact
groupoids, we introduce the following definition.
\begin{definition}\label{definicion}
Let $\gr$ be a Lie groupoid, $(\Jacobi )$ be a Jacobi structure on
$G$ and $\sigma :G\to \R$ be a multiplicative function. Then,
$(\gr,\Jacobi ,\sigma )$ is a {\em Jacobi groupoid} if the
homomorphism $\sostJ:T^\ast G\times \R\to TG\times \R$ given by
$$\sostJ (\omega _g,\gamma )=(\sostP(\omega _g)+\gamma
\,E_g,-\omega _g(E_g))$$is a morphism of Lie groupoids over some
map $\varphi _0:A^\ast G\to TM\times \R$, where the structural
functions of the Lie groupoid structure on $T^\ast G\times \R\gpd
A^\ast G$ (respectively, $TG\times \R \gpd TM\times \R$) are given
by (\ref{T*GR}) (respectively, (\ref{TGR})).
\end{definition}
\begin{remark}
{\rm Since $\sostJ:T^\ast G\times \R\to TG\times \R$ is a morphism
of Lie groupoids, we deduce that $$\varphi _0=(\s ^T)_\sigma \circ
\sostJ\circ \tilde{\epsilon}_\sigma =(\t ^T)_\sigma \circ
\sostJ\circ \tilde{\epsilon}_\sigma .$$ Thus, if $\omega _x\in
A^\ast _xG$, it follows that
\begin{equation}\label{aplbase}
\varphi _0(\omega _x)=\Big ( \s ^{\tilde{x}}_\ast ( \sostP
(\tilde{\epsilon}(\omega _x))),-\omega _x\,(E(\tilde{x})-\epsilon
_\ast ^x(\t _\ast ^{\tilde{x}} (E(\tilde{x}))))\Big ).
\end{equation}
}
\end{remark}
\begin{examples}\label{ejemplos}
{\rm {\bf 1.-}{\it Poisson groupoids}

  If $(\gr,\Jacobi ,\sigma
)$ is a Jacobi groupoid with $E=0$ and $\sigma =0$ then we recover
the definition of a Poisson groupoid (see \cite{MX,MX2} and Remark
\ref{recup}).

{\bf 2.-}{\it Contact groupoids}

 Let $(\gr ,\eta ,\sigma )$ be a
contact groupoid. If $(\Jacobi )$ is the Jacobi structure
associated with the contact 1-form $\eta$ then, using the results
in Section \ref{caract}, we have that $(\gr, \Jacobi ,\sigma )$ is
a Jacobi groupoid.

{\bf 3.-}{\it Jacobi-Lie groups}

In \cite{IM2}, we proved that generalized Lie bialgebras (that is,
generalized Lie bialgebroids over a single point) may be
considered as the infinitesimal invariants of a particular class
of Lie groups. These Lie groups can be defined as follows. Let $G$
be a Lie group with identity element $\mathfrak e$, $\mult$ be a
multiplicative function and $(\Jacobi )$ be a Jacobi structure on
$G$ such that:
\begin{enumerate}
\item $\Lambda$ is $\sigma$-multiplicative, i.e.,
$\Lambda (gh)=(R_h)^g_\ast (\Lambda (g))+e^{-\sigma
(g)}(L_g)^h_\ast (\Lambda (h)),$ for $g,h\in G$.
\item $E$ is a right-invariant vector field, $E(\mathfrak
e)=-X_0$.
\item $\sostP (\delta \sigma )=\der{X}_0-
e^{-\sigma}\izq{X}_0.$
\end{enumerate}
Condition (i) implies that $\Lambda (\mathfrak e )=0$ and
conditions (ii) and (iii) imply that $E(\sigma )=0$. Thus, using
again (ii) and (iii), we deduce that $$(\s ^T)_\sigma \circ \sostJ
=\varphi _0\circ \tilde{\s}_\sigma,\quad (\t ^T)_\sigma \circ
\sostJ =\varphi _0\circ \tilde{\t}_\sigma .$$ In addition, from
condition (ii) and (iii), we have that
\begin{equation}\label{laE}
(L_g)_\ast ^hE_h=e^{\sigma (g)}( E_{gh}+(R_h)_\ast ^g(\sostP
(\delta \sigma )_g)),
\end{equation}
for $g,h\in G$.

Now, suppose that $(\omega _g,\gamma )\in T^\ast_g G\times\R$ and
$(\nu _h,\zeta )\in T^\ast _hG\times\R$ satisfy the condition
$\tilde{\s}_\sigma (\omega _g,\gamma )=\tilde{\t}_\sigma (\nu
_h,\zeta )$. Then,
\begin{equation}\label{compati}
\tilde{\s}(\omega _g+\zeta\, e^{\sigma (g)}(\delta \sigma
)_g)=\tilde{\t}(e^{\sigma (g)}\nu _h).
\end{equation}
Thus, using (\ref{mult'}), (\ref{TGR}), (\ref{laE}) and the fact
that $E$ is a right-invariant vector field, we deduce that
$$\begin{array}{l} \sostJ (\omega _g,\gamma )\oplus _{TG\times\Rp}
\sostJ(\nu _h,\zeta )\kern2pt =\\[5pt] \kern20pt \Big ((L_g)_\ast
^h(\sostP (\nu _h))+(R_h)_\ast ^g(\sostP(\omega _g+\zeta \,
e^{\sigma (g)}(\delta \sigma )_g)+(\gamma +e^{\sigma (g)}\zeta
)E_{gh},-\omega _g(E_g)\Big ).\end{array}$$

On the other hand, $E_{gh}=E_g\oplus _{TG}0_h,$ and therefore,
from (\ref{mult*}) and (\ref{T*GR}), it follows that
$$\begin{array}{l}\sostJ ((\omega _g,\gamma )\oplus _{T^\ast
G\times \Rp}(\nu _h,\zeta ))\kern2pt =\\[5pt] \kern20pt
\displaystyle \Big ( \frac{1}{2}\sostP \Big \{ ((R_{h\inv})_\ast
^{gh})^\ast (\omega _g+\zeta \,e^{\sigma (g)}(\delta \sigma )_g)+
((L_{g\inv})_\ast ^{gh})^\ast (e^{\sigma (g)}\nu _h)\Big \}
\\[5pt] \kern25pt +(\gamma +e^{\sigma (g)}\zeta )E_{gh},-\omega
_g(E_g)\Big ).\end{array}$$

Consequently, using (\ref{T*G}), (\ref{compati}) and the fact that
$\Lambda$ is $\sigma$-multiplicative, we conclude that $$\sostJ
(\omega _g,\gamma )\oplus _{TG\times\Rp} \sostJ(\nu _h,\zeta
)=\sostJ ((\omega _g,\gamma )\oplus _{T^\ast G\times \Rp}(\nu
_h,\zeta )).$$ Thus, we have proved that $(\gr,\Lambda ,E,\sigma
)$ is a Jacobi groupoid.

{\bf 4.-}{\it An abelian Jacobi groupoid}

Let $(L,\lcf \, ,\, \rcf ,\rho )$ be a Lie algebroid over $M$ and
$\Lambda_{L^\ast }$ be the corresponding linear Poisson structure
on the dual bundle $L^\ast$ (see Section \ref{algebroides}). We
may consider on $L^\ast$ the Lie groupoid structure for which $\s
=\t$ is the vector bundle projection and the partial
multiplication is the addition in the fibers. Then, $L^\ast$ with
the Poisson structure $\Lambda _{L^\ast}$ is a Poisson groupoid
(see \cite{We2}).

Now, suppose that $\omega _0\in \Gamma (L^\ast)$ is a 1-cocycle of
$L$ and denote by $(\Lambda_{(L^\ast ,\omega _0)},E_{(L^\ast
,\omega _0)})$ the Jacobi structure on $L^*$ given by
(\ref{Jacobi-lineal}). Note that: i) The Liouville vector field
$\Delta$ of $L^\ast$ and the vertical lift $\omega _0^v\in
\mathfrak X (L^\ast)$  of $\omega _0$ to $L^\ast$ are
$\s$-vertical and $\t$-vertical vector fields on $L^\ast,$ and ii)
$\omega_0^v$ is a right-invariant and left-invariant vector field
on $L^\ast$. Using i), ii), (\ref{sumaTG}), (\ref{TGR}),
(\ref{T*GR}), (\ref{Jacobi-lineal}) and the fact that $(L^\ast
,\Lambda _{L^\ast})$ is a Poisson groupoid, we deduce that
$(L^\ast \gpd M,\Lambda_{(L^\ast ,\omega _0)},E_{(L^\ast ,\omega
_0)},0)$ is a Jacobi groupoid.

{\bf 5.-}{\it The banal Jacobi groupoid}

  Let $M$ be a
differentiable manifold. The results in Section \ref{grupoides}
(see Examples \ref{ej-grupoides}) imply that $G= M\times \R \times
M$ is a Lie groupoid over $M$ and, moreover, the function $\mult$
given by $\sigma (x,t,y)=t$ is multiplicative. Thus, we can
consider the corresponding Lie groupoids $TG\times \R\gpd
TM\times\R$ and $T^\ast G\times\R\gpd A^\ast G$.

On the other hand, the map $\Phi :TM\times\R \to AG$ given by
\begin{equation}\label{2.9'}
\Phi (X_x,\lambda )=(0,\lambda \frac{\partial }{\partial t}_{|
0},X_x)\in T_{(x,0,x)}G,\mbox{ for }(X_x,\lambda )\in
T_xM\times\R,
\end{equation}
defines an isomorphism between the Lie algebroids $(TM\times \R
,\makebox{{\bf [}\, ,\, {\bf ]}}, \pi)$ (see Section
\ref{algebroides}) and $AG$. Thus, $AG$ may be identified with
$TM\times\R$ and, under this identification, the projections and
the partial multiplications on $TG\times\R$ and $T^\ast G\times\R$
are given by $$\begin{array}{l} (\s ^T)_\sigma ((X_x,a\fract
_{|t},Y_y),\lambda )=(Y_y,a+\lambda ),\\ (\t ^T)_\sigma
((X'_{x'},a'\fract _{|t'},Y'_{y'}),\lambda ')=(X'_{x'},\lambda
'),\\ ((X_x,a\fract _{|t},Y_y),\lambda
)\oplus_{TG\times\Rp}((Y_y,a'\fract _{|t'},Y'_{y'}),a+\lambda
)=((X_x,(a+a')\fract _{|t+t'},Y'_{y'}),\lambda ),\\[5pt]
\tilde{\s} _\sigma ((\omega _x,a\,\delta t _{|t},\theta _y),\gamma
)=(e^{-t}\theta _y,\gamma ),\\ \tilde{\t} _\sigma ((\omega
'_{x'},a'\,\delta t_{|t'},\theta '_{y'}),\gamma ')=(-\omega
'_{x'},a'-\gamma '),\\((\omega _x,a\,\delta t _{|t},\theta
_y),\gamma )\oplus _{T^\ast G\times\Rp} ((-e^{-t}\theta
_{y},a'\,\delta t_{|t'},\theta
'_{y'}),a'-e^{-t}a)\kern2pt=\\\kern220pt((\omega _x,a'e^t\,\delta
t_{|t+t'},e^t\theta '_{y'}),\gamma -a+e^t\,a').
\end{array}$$
Now, suppose that $(\Jacobi )$ is a Jacobi structure on $M$. Then,
it was proved in \cite{ILMM} that the pair $(\Lambda ',E')$ is a
Jacobi structure on $G$, where
\begin{equation}\label{2.9''}
\Lambda '(x,t,y)=-\Big (\Lambda (x)-\fract _{|t}\wedge E(x)\Big
)+e^{-t}\Big ( \Lambda (y)+\fract _{|t}\wedge E(y)\Big ),\quad
E'(x,t,y)=-E(x).
\end{equation}
Furthermore, it is easy to prove that the map $\varphi _0:A^\ast
G\cong T^\ast M\times\R\to TM\times\R$ given by (\ref{aplbase}) is
just the homomorphism $\sostJ :T^\ast M\times\R\to TM\times\R$.
Using the above facts, we conclude that $(\gr , \Lambda
',E',\sigma )$ is a Jacobi groupoid.
 }
\end{examples}
\subsection{Some basic properties of Jacobi groupoids}
In this Section, we will show some basic properties of Jacobi
groupoids.
\begin{proposition}\label{propiedades}
Let $(\gr ,\Jacobi ,\sigma )$ be a Jacobi groupoid. Then:
\begin{itemize}
\item[{\it i)}] $M\cong \epsilon (M)$ is a coisotropic submanifold
in $G$.
\item[{\it ii)}] $E$ is a right-invariant vector field on
$G$ and $E(\sigma )=0$. Moreover, if $X_0\in \Gamma (AG)$ is the
section of the Lie algebroid $AG$ of $G$ satisfying
$E=-\der{X_0}$, we have that
\begin{equation}\label{compatibilidad}
\sostP (\delta \sigma )= \der{X_0}-e^{-\sigma}\izq{X_0}.
\end{equation}

\vspace{-.5cm}

\item[{\it iii)}] If $\tilde{\s}$, $\tilde{\t}$ and $\tilde{\epsilon}$
are the projections and the inclusion of the Lie groupoid $T^\ast
G\gpd A^\ast G$ then $$e^{-\sigma}\sostP\circ
\tilde{\epsilon}\circ\tilde{\s}=\epsilon ^T\circ \s ^T\circ
\sostP,\quad \sostP\circ\tilde{\epsilon}\circ\tilde{\t}=\epsilon
^T\circ \t ^T\circ \sostP.$$

\vspace{-.25cm}
\item[{\it iv)}] If $g$ and $h$ are elements of $G$ such that
$\s (g)= \t (h)=x$ and ${\cal X}$ and ${\cal Y}$ are (local)
bisections through the points $g$ and $h$, ${\cal X}(x)=g$ and
${\cal Y}(x)=h$, then
\begin{equation}\label{sigma-multi}
\Lambda (gh) = (R_{{\cal Y} })_\ast ^g (\Lambda (g)) + \expsn
(L_{{\cal X}})_\ast ^h (\Lambda (h)) - \expsn (L_{{\cal X}} \circ
R_{{\cal Y}})_\ast ^{\tilde{x}} (\Lambda (\tilde{x} )).
\end{equation}
\end{itemize}
\end{proposition}
\prueba If $x$ is a point of $M$ then,  using  (\ref{T*G}), we
obtain that the map $$\tilde{\epsilon}_{|A^\ast _xG}:A^\ast _xG\to
T^\ast _{\tilde{x}}G$$ is a linear isomorphism between the vector
spaces $A^ \ast _xG$ and the annihilator of the subspace
$T_{\tilde{x}}\epsilon (M)$, that is, $(T_{\tilde{x}}\epsilon
(M))^\circ$. Thus, from (\ref{TG}), (\ref{T*G}), (\ref{TGR}),
(\ref{T*GR}), (\ref{aplbase}) and since $(\epsilon^T)_\sigma \circ
\varphi_0=\#_{(\Jacobi)}\circ \tilde{\epsilon}_\sigma,$ it follows
that $M\cong \epsilon (M)$ is a coisotropic submanifold in $G$
with respect to $\Lambda$. This proves {\it i)}.

On the other hand, using (\ref{TG}), (\ref{T*G}), (\ref{TGR}),
(\ref{T*GR}), (\ref{aplbase}) and the relations $$(\s ^T)_\sigma
\circ \sostJ =\varphi _0\circ \tilde{\s}_\sigma,\quad (\t
^T)_\sigma \circ \sostJ =\varphi _0\circ \tilde{\t}_\sigma ,$$ we
deduce {\it ii)} and {\it iii)}.

Finally, we will prove {\it iv)}. Using the multiplicative
function $\sigma$, one may introduce the Lie groupoid structure in
$T^\ast G$ over $A^\ast G$ with structural functions
$\tilde{\s}^\ast _\sigma$, $\tilde{\t}_\sigma ^\ast$, $\oplus
^\sigma _{T^\ast G}$, $\tilde{\epsilon}_\sigma ^\ast$ and
$\tilde{\iota} ^\ast _\sigma$ given by
\begin{equation}\label{T*Gsigma}
\begin{array}{l}
\tilde{\s}^\ast _\sigma(\omega _g
)=e^{-\sigma(g)}\tilde{\s}(\omega _g), \quad \tilde{\t}^\ast
_\sigma(\nu _h)=\tilde{\t}(\nu _h),\mbox{ for }\omega _g\in T^\ast
_gG\mbox{ and }\nu _h\in T^\ast _hG,\\[5pt] ( \omega _g \oplus
^\sigma _{T^\ast G}\nu _h )= \omega _g \oplus _{T^\ast
G}(e^{\sigma (g)}\nu _h),\\[5pt] \tilde{\epsilon}_\sigma ^\ast
(\omega _x)=\tilde{\epsilon} (\omega _x),\mbox{ for }\omega _x\in
A^\ast _x G,\\[5pt] \tilde{\iota} ^\ast _\sigma (\omega
_g)=e^{-\sigma (g)}\tilde{\iota}(\omega _g),\mbox{ for }\omega
_g\in T^\ast _gG.
\end{array}
\end{equation}
In fact, if we consider on $T^\ast G\times \R$ the Lie groupoid
structure over $A^\ast G$ introduced in Section \ref{caract} then
the canonical inclusion $$T^\ast G\to T^\ast G\times \R,\quad
\omega _g\in T^\ast _gG\mapsto (\omega _g,0)\in T^\ast
_gG\times\R,$$ is a Lie groupoid monomorphism over the identity of
$A^\ast G$.

Since the map $\sostJ: T^\ast G\times\R\to TG\times\R$ is a Lie
groupoid homomorphism, we have that (see (\ref{TGR}), (\ref{T*GR})
and (\ref{T*Gsigma})) $$\sostP ( \omega _g \oplus_{T^\ast G}\nu _h
)=\sostP (\omega _g)\oplus _{TG}\sostP (e^{-\sigma(g)}\nu _h)$$
for $\omega _g\in T^\ast _gG$ and $\nu _h\in T^\ast _hG$
satisfying $\tilde{\s}(\omega _g )=\tilde{\t}(\nu _h)$. Thus, if
$\Pi$ is the 2-vector on $G\times G\times G$ defined by $\Pi
(g,h,k)=e^{\sigma (g)}\Lambda (g)+\Lambda (h)-e^{\sigma
(g)}\Lambda (k),$ it follows that the graph of the multiplication
in $G,$ $\{ (g,h,gh)\in G\times G\times G\, /\, \s (g)=\t(h)\},$
is a coisotropic submanifold of $G\times G\times G$ with respect
to $\Pi$.

Now, denote by $\Omega$ the affinoid diagram corresponding to the
Lie groupoid $G$, that is (see \cite{We3}), $$\Omega = \{
(k,g,h,r)\in G\times G\times G\times G \, /\, \s (h)=\s (k),\, \t
(k)=\t (g),\, r=hk\inv g \} .$$ Then, following the proof of
Theorem 4.5 in \cite{We3}, we obtain that $\Omega$ is a
coisotropic submanifold of $G\times G\times G\times G$ with
respect to the 2-vector $\widetilde{\Pi}$ given by
$$\widetilde{\Pi} (k,g,h,r)=e^{\sigma (k)}\Lambda (k)-e^{\sigma
(k)}\Lambda (g)-e^{\sigma(h)}\Lambda (h)+e^{\sigma (h)}\Lambda
(r).$$ On the other hand, if $g$ and $h$ are elements of $G$
satisfying $\s (g)=\t (h)=x$, we have that $(gh, g,h,\tilde{x})$
is an element of $\Omega$. In addition, for any $\xi \in T^\ast
_{gh}G$ and ${\cal X},{\cal Y}$ (local) bisections of $G$ through
the points $g$ and $h$ (${\cal X}(x)=g$ and ${\cal Y}(x)=h$), it
follows from Lemma 2.6 in \cite{Xu} that $$(-\xi ,((R_{\cal
Y})_\ast ^g)^\ast (\xi),((L_{\cal X})_\ast ^h)^\ast (\xi
),-((R_{\cal Y}\circ L_{\cal X})_\ast ^{\tilde{x}})^\ast (\xi ))$$
is a conormal vector to $\Omega$ at $(gh,g,h,\tilde{x})$, i.e., it
is an element of $(T_{(gh,g,h,\tilde{x})}\Omega )^\circ$.
Therefore, if $\xi,\eta \in T^\ast _{gh}G$, we deduce that $$\Big
( e^{\sigma (gh)}\Lambda (gh)-e^{\sigma (h)}(L_{\cal X})_\ast
^h(\Lambda (h))-e^{\sigma (gh)}(R_{\cal Y})_\ast ^g(\Lambda
(g))+e^{\sigma (h)} (R_{\cal Y}\circ L_{\cal X})_\ast
^{\tilde{x}}(\Lambda (\tilde{x}))\Big ) (\xi ,\eta )=0.$$ This
implies that (\ref{sigma-multi}) holds.\QED

Motivated by the above result, we introduce the following
definition.
\begin{definition}
Let $\gr$ be a Lie groupoid and $\mult$ be a multiplicative
function. A multivector field $P$ on $G$ is {\em $\sigma$-affine}
if for any $g, h\in G$ such that $\s (g)= \t (h)=x$ and any
(local) bisections ${\cal X}, {\cal Y}$ through the points $g, h$,
${\cal X}(x)=g$ and ${\cal Y}(x)=h$, we have
\begin{equation}\label{eq:multi}
P(gh) = (R_{{\cal Y} })_\ast ^g (P(g)) + \expsn (L_{{\cal
X}})_\ast ^h (P (h)) - \expsn (L_{{\cal X}}\circ R_{{\cal
Y}})_\ast ^{\tilde{x}} (P (\tilde{x}) ).
\end{equation}
\end{definition}
It is clear that if $P$ is a $\sigma$-affine multivector and
$\sigma$ identically vanishes, then $P$ is affine (see
\cite{MX2,Xu}). On the other hand, if $G$ is a Lie group with
identity element $\mathfrak e$ and $P$ is a $\sigma$-affine
multivector field on $G$ such that $P(\mathfrak e)=0$, then $P$ is
a $\sigma$-multiplicative multivector field in the sense of
\cite{IM2}.

The following proposition gives a very useful characterization of
$\sigma$-affine multivector fields (see \cite{IM2} for the
corresponding result for the case of Lie groups).
\begin{proposition}\label{sigma-afin}
Let $\gr$ be an $\s$-connected Lie groupoid and $\mult$ be a
multiplicative function on $G$. For a multivector field $P$ on
$G$, the following statements are equivalent:
\begin{enumerate}
\item $P$ is $\sigma$-affine;
\item For any left-invariant vector field $\izq{X}$, the Lie derivative
$e^\sigma {\cal L}_{\ssli} P$ is left-invariant.
\end{enumerate}
\end{proposition}
\prueba The result follows using the fact that $\sigma$ is
multiplicative and proceeding as in the proof of Theorem 2.2 in
\cite{MX2}.\QED
\section{Jacobi groupoids and generalized Lie bialgebroids}
\setcounter{equation}{0} The aim of this Section is to show the
relation between Jacobi groupoids and generalized Lie
bialgebroids.
\subsection{Coisotropic submanifolds of a Jacobi manifold, Lie
algebroids and 1-cocycles} In this Section, we will prove that if
$S$ is a coisotropic submanifold of a Jacobi manifold $M$ then
there exists a Lie algebroid structure on the conormal bundle to
$S$ and, in addition, we can define a distinguished 1-cocycle for
this Lie algebroid structure. For this purpose, we will need the
following result.

\begin{lemma}\label{lema-coisotrop}
Let $(M, \Jacobi )$ be a Jacobi manifold and $(\lcf \, ,\,\rcf
_{(\Jacobi )},$ $\widetilde{\#}_{(\Jacobi )})$ be the Lie
algebroid structure  on $T^\ast M \times \R$. Suppose that $S$ is
a coisotropic submanifold of $M$ and that $\bar{\jmath
}^\ast:\Omega ^1(M)\times C^\infty (M,\R)\to \Omega ^1(S)\times
C^\infty (S,\R)$ is the map defined by $\bar{\jmath }^\ast(\omega
,f)=(\jmath ^\ast \omega ,\jmath^* f)$, $\jmath :S\to M$ being the
canonical inclusion. Then:
\begin{itemize}
\item[{\it i)}] $Ker\, \bar{\jmath }^\ast$ is a Lie subalgebra of
the Lie algebra $(\Omega ^1(M)\times C^\infty (M,\R),\lcf\, ,\,
\rcf _{(\Jacobi)})$.
\item[{\it ii)}] The subspace of $\,\Omega ^1(M)\times C^\infty (M,\R)$
defined by $\{ (\omega ,f)\in \Omega ^1(M)\times C^\infty (M,\R)
\, / \, \omega _{|S}=0, \jmath^* f=0\}$ is an ideal in $Ker\,
\bar{\jmath }^\ast$.
\end{itemize}
\end{lemma}
\prueba {\it i)} If $(\omega ,f)$, $(\nu ,g )\in \Omega
^1(M)\times C^\infty (M,\R )$ satisfy $$\bar{\jmath }^\ast (\omega
,f)=0,\quad \bar{\jmath }^\ast(\nu ,g )=0,$$ it follows from
(\ref{ecjacobi}) that
\begin{equation}\label{0-ec}
\begin{array}{lll}
\bar{\jmath }^\ast \lcf (\omega ,f), (\nu ,g )\rcf _{(\Jacobi
)}&=&(\jmath ^\ast (i(\sostP (\omega ))\delta \nu - i(\sostP (\nu
))\delta \omega -\delta (\omega (\sostP (\nu )))),\\
&&\jmath^*(\omega (\sostP (\nu ))+ \sostP(\omega )(g)-\sostP (\nu
)(f))).
\end{array}
\end{equation}
Now, since $\jmath^*\omega=0, \jmath^*\nu=0$ and $S$  is a
coisotropic submanifold,  it follows that the restriction to $S$
of the vector fields $\sostP (\omega )$ and $\sostP (\nu )$ is
tangent to $S$. Thus, from (\ref{0-ec}), we deduce that
$$\bar{\jmath }^\ast \lcf (\omega ,f), (\nu ,g )\rcf _{(\Jacobi
)}=0.$$ {\it ii)} If $\omega '$ and $\nu '$ are 1-form on $M$, we
will denote by $\lcf \omega ',\nu '\rcf _{\Lambda}$ the 1-form on
$M$ given by $$\lcf \omega ',\nu '\rcf _{\Lambda}=i(\sostP (\omega
' ))\delta \nu '- i(\sostP (\nu '))\delta \omega '-\delta (\omega
' (\sostP (\nu '))).$$ Note that
\begin{equation}\label{1-ec}
\lcf \omega ',f\nu '\rcf _{\Lambda}=f\lcf \omega ',\nu '\rcf
_{\Lambda}+\sostP(\omega ')(f)\nu ',\mbox{ for }f\in C^\infty
(M,\R).
\end{equation}
Next, suppose that $(\omega ,f),(\nu ,g)\in \Omega ^1(M)\times
C^\infty (M,\R )$ satisfy the following conditions $$\omega
_{|S}=0,\quad \jmath^* f=0,\quad \bar{\jmath}^\ast (\nu ,g)=0.$$
Then, proceeding as in the proof of {\it i)}, we have that $$\lcf
(\omega ,f), (\nu ,g )\rcf _{(\Jacobi )}{}_{|S}=(\lcf \omega ,\nu
\rcf _{\Lambda}{}_{|S},0).$$ Thus, if $x$ is a point of $S$, we
must prove that $\lcf \omega ,\nu \rcf _{\Lambda}(x)=0$. For this
purpose, we consider a coordinate neighborhood $(U,\varphi )$ of
$M$ with coordinates $(x_1,\ldots ,x_n,x_{n+1},\ldots ,x_m)$ such
that $$\varphi (U\cap S)=\{ (x_1,\ldots ,x_m) \in \varphi (U)\,
/\, x_{n+1}=\ldots =x_m=0\} .$$ Here, $n$ (respectively, $m$) is
the dimension of $S$ (respectively, $M$). Then, on $U$
\begin{equation}\label{2-ec}
\omega = \sum _{i=1}^m \omega ^i\,\delta x_i,\qquad \nu = \sum
_{j=1}^n \nu ^j\,\delta x_j+\sum _{k=n+1}^m\bar{\nu }^k\,\delta
x_k
\end{equation}
with
\begin{equation}\label{3-ec}
\jmath^*\omega ^i=0,\quad \jmath^*\nu ^j=0,
\end{equation}
for all $i\in \{ 1,\ldots ,m\}$ and $j\in \{ 1,\ldots ,n\}$.

Note that, since $S$ is a coisotropic submanifold of $M$, it
follows that
\begin{equation}\label{4-ec}
\sostP(\delta x_k )_{|S}(\omega ^i)=0,\mbox{ for all }i\in \{
1,\ldots ,m\}\mbox{ and }k\in \{ n+1,\ldots ,m\}.
\end{equation}
Therefore, using (\ref{1-ec})-(\ref{4-ec}), we conclude that $\lcf
\omega ,\nu \rcf _{\Lambda}(x)=0$.\QED

Now, we will show the main result of the Section.
\begin{proposition}\label{alg-ind}
Let $(M, \Jacobi )$ be a Jacobi manifold and $(\lcf \, ,\,\rcf
_{(\Jacobi )},$ $\widetilde{\#}_{(\Jacobi )})$ be the Lie
algebroid structure  on $T^\ast M \times \R$. Suppose that $S$ is
a coisotropic submanifold of $M$. Then:
\begin{itemize}
\item[{\it i)}] The conormal bundle to $S$, $N(S)=(TS)^\circ \to S$,
admits a Lie algebroid structure $(\lcf \, ,\, \rcf _S,\rho _S)$
defined by
\begin{equation}\label{corch-coi}
\begin{array}{lll} \lcf \omega ,\nu \rcf _S(x) &=&(\pi _1 \lcf
(\tilde{\omega},0) ,(\tilde{\nu},0) \rcf _{(\Jacobi)})(x),\\ \rho
_S(\omega )(x)&=&\sostP(\omega _x),
\end{array}
\end{equation}
for all $x\in S$, where $\pi _1:\Omega ^1(M)\times C^\infty (M,\R
)\to \Omega^1 (M)$ is the projection onto the first factor and
$\tilde{\omega}$ and $\tilde{\nu}$ are arbitrary extensions to $M$
of $\,\omega$ and $\nu$, respectively.
\item[{\it ii)}] The section $E_S$ of the vector bundle $N(S)^\ast\to S$
characterized by
\begin{equation}\label{1coci-coi}
\omega (E_S(x))=-\omega (E(x))
\end{equation}
for all $\omega \in N_xS=(T_xS)^\circ$ and $x\in S$, is a
1-cocycle of the Lie algebroid $(N(S),\lcf \, ,\, \rcf _S,$ $\rho
_S)$.
\end{itemize}
\end{proposition}
\prueba  {\it i)} follows from Lemma \ref{lema-coisotrop} and {\it
ii)} follows using (\ref{1coci-coi}) and the fact that $(-E,0)\in
\mathfrak X (M)\times C^\infty (M,\R )$ is a 1-cocycle of the Lie
algebroid $(T^\ast M \times \R,\lcf \, ,\,\rcf _{(\Lambda
,E)},\widetilde{\#}_{(\Lambda ,E)})$.\QED
\begin{remark}\label{cocy0}
{\rm If the Jacobi manifold $M$ is Poisson (that is, $E=0$) then
the 1-cocycle $E_S$ identically vanishes and $(\lcf \, ,\, \rcf
_S,\rho _S)$ is just the Lie algebroid structure obtained by
Weinstein in \cite{We2}.}
\end{remark}
\subsection{The generalized Lie bialgebroid of a Jacobi groupoid}
In this Section, we will show that generalized Lie bialgebroids
are the infinitesimal invariants for Jacobi groupoids.

Let $(G\gpd M,\Jacobi ,\sigma )$ be a Jacobi groupoid and $AG$ be
the Lie algebroid of $G$. Then, $E$ is a right-invariant vector
field and, thus, there exists a section $X_0$ of $AG$ such that
$E=-\der{X_0}$ (see Proposition \ref{propiedades}). Moreover, the
conormal bundle to $M$, as a submanifold of $G$, may be identified
with $A^\ast G$.

Now, we consider the section $\phi _0$ of $A^\ast G$ given by
\begin{equation}\label{fi0}
\phi _0(X_x)=X _x(\sigma ), \mbox{ for }X_x\in A_xG\mbox{ and
}x\in M.
\end{equation}
Since $\sigma$ is a Lie groupoid 1-cocycle, it follows that $\phi
_0$ is a 1-cocycle of the Lie algebroid $AG$ (see \cite{Xu}).

On the other hand, using that $M\cong \epsilon (M)$ is a
coisotropic submanifold of $G$, we deduce that there exists a Lie
algebroid structure $(\lcf \, ,\, \rcf _\ast ,\rho _\ast )$ on
$A^\ast G$ and, furthermore, the vector field $E$ induces a
1-cocycle $E_M\in \Gamma (AG)$ of $A^\ast G$ (see Proposition
\ref{alg-ind}). In fact, from Proposition \ref{alg-ind}, we have
that $E_M=X_0$ and
\begin{equation}\label{3.8I}
\begin{array}{rcl}
\lcf \omega,\nu\rcf_\ast(x)&=&\pi_1\lcf
(\widetilde{\tilde\epsilon\circ
\omega},0),(\widetilde{\tilde\epsilon\circ
\nu},0)\rcf_{(\Jacobi)}(\tilde{x}),\;\;\; \rho_*(\omega)(x)=\s
_*^{\tilde{x}} (\#_\Lambda(\tilde{\epsilon}(\omega_x))),\\
\end{array}\end{equation}
for $\omega,\nu\in \Gamma(A^*G)$ and $x\in M$, where
$\tilde{\epsilon}$ is the inclusion in the Lie groupoid $T^\ast
G\gpd A^\ast G$ and $\widetilde{\tilde\epsilon\circ \omega}$ and
$\widetilde{\tilde\epsilon\circ \nu}$ are arbitrary extensions to
$G$ of $\tilde\epsilon\circ \omega$ and $\tilde\epsilon\circ \nu$,
respectively.

Note that, from (\ref{aplbase}) and (\ref{3.8I}), we have that
$\varphi _0=(\rho _\ast ,X_0)$, where
\begin{equation}\label{3.8II}
(\rho _\ast ,X_0)(\omega _x)=(\rho _\ast (\omega _x),\omega _x(X_0
(x)))
\end{equation}
for $\omega _x\in \displaystyle A_x^\ast G$.
 In addition, we will prove the following
result.
\begin{theorem}\label{bajada}
Let $(G\gpd M,\Jacobi ,\sigma )$ be a Jacobi groupoid. Then
$((AG,\phi _0),(A^\ast G,X_0))$ is a generalized Lie bialgebroid.
\end{theorem}
\prueba Denote by $d_{\ast X_0}$ the $X_0$-differential of the Lie
algebroid $(A^\ast G,\lcf \, ,\, \rcf _\ast ,\rho _\ast )$.

We will show that
\begin{equation}\label{dif-induc}
e^\sigma {\cal L}_\ssli \Lambda =-\izq{d_{\ast X_0}X}.
\end{equation}
for $X\in \Gamma (AG)$. Suppose that $\omega _1,\omega _2$ are any
sections of $A^\ast G$. Let $\widetilde{\tilde\epsilon\circ
\omega} _1,\widetilde{\tilde\epsilon\circ \omega}_2$ be any of
their extensions to 1-forms on $G$. Then, using
(\ref{diferencial}), (\ref{ecjacobi}), (\ref{ndif}), (\ref{3.8I})
and the fact that $\sigma _{|\epsilon (M)}\equiv 0$, we have that
$$\begin{array}{lll} \Big ( e^\sigma {\cal L}_\ssli \Lambda\Big
)_{|\epsilon (M)}(\omega _1,\omega _2 )&=& \Big ( ({\cal
L}_{\sostP (\widetilde{\tilde\epsilon\circ \omega}_1)}
\widetilde{\tilde\epsilon\circ \omega} _2-{\cal L}_{\sostP
(\widetilde{\tilde\epsilon\circ \omega} _2)}
\widetilde{\tilde\epsilon\circ \omega}_1 - \Lambda
(\widetilde{\tilde\epsilon\circ \omega}
_1,\widetilde{\tilde\epsilon\circ \omega} _2))(\izq{X})\\&&+\sostP
(\widetilde{\tilde\epsilon\circ \omega} _2 )
(\widetilde{\tilde\epsilon\circ \omega} _1(\izq{X}))-\sostP
(\widetilde{\tilde\epsilon\circ \omega} _1
)(\widetilde{\tilde\epsilon\circ \omega}_2(\izq{X})) \Big
)_{|\epsilon (M)}\\ &=& \lcf\omega _1,\omega _2\rcf _\ast (X)+\rho
_\ast (\omega _2)(\omega _1(X))\\&&-\rho _\ast (\omega _1)(\omega
_2(X))-(X_0\wedge X)(\omega _1, \omega _2)\\&=& -(d_{\ast
X_0}X)(\omega _1, \omega _2).
\end{array}$$
Thus, since $-\izq{d_{\ast X_0}X}$ and $e^\sigma {\cal L}_\ssli
\Lambda$ are left-invariant 2-vectors (see Proposition
\ref{sigma-afin}) and their evaluation coincides on the conormal
bundle $A^\ast G,$ we deduce (\ref{dif-induc}).

Using (\ref{multi-inv-izq}), (\ref{fi0}) and (\ref{dif-induc}), we
obtain that
\begin{equation}\label{ecu1}
\begin{array}{lll}
\izq{d_\ast{}_{X_0}\lcf X ,Y \rcf }&=&-e^\sigma {\cal L}_{[\ssli,
\sslib ]}\Lambda\\&=& {\cal L}_\sslib (e^\sigma {\cal L}_\ssli
\Lambda )- \izq{Y}(\sigma )(e^\sigma {\cal L}_\ssli \Lambda )
-{\cal L}_\ssli (e^\sigma {\cal L}_\sslib \Lambda )+
\izq{X}(\sigma )(e^\sigma {\cal L}_\sslib \Lambda )\\&=&\izq{\lcf
X ,d_\ast{}_{X_0}Y \rcf} - \izq{\phi _0(X)d_\ast{}_{X_0}Y}
-\izq{\lcf Y ,d_\ast{}_{X_0}X \rcf}+\izq{\phi
_0(Y)d_\ast{}_{X_0}X},
\end{array}
\end{equation}
for $X,Y\in\Gamma (AG)$, where $\estalg$ is the Lie algebroid
structure on $AG$. Thus, from (\ref{Schouext}) and (\ref{ecu1}),
we conclude that $$d_\ast{}_{X_0}\lcf X ,Y \rcf = \lcf X
,d_\ast{}_{X_0}Y \rcf _{\phi _0}-\lcf Y ,d_\ast{}_{X_0}X \rcf
_{\phi _0}$$ for $X,Y \in \Gamma (AG)$.

Now, (\ref{fi0}), the condition $E(\sigma )=-\der{X_0}(\sigma )=0$
(see Proposition \ref{propiedades}) and the fact that $\sigma$ is
a multiplicative function imply that $\phi _0(X_0)\circ \s =0$
and, therefore,
\begin{equation}\label{util}
\phi _0(X_0)=0.
\end{equation}
Furthermore, if $x\in M$ then, from (\ref{compatibilidad}),
(\ref{fi0}) and (\ref{3.8I}), we deduce that $$\epsilon^x_\ast
(\rho_*(\phi_0)(x))=\sostP(\delta \sigma
)(\tilde{x})=\izq{X_0}(\tilde{x})-\der{X_0}(\tilde{x})=-\epsilon
^x_\ast (\s ^{\tilde{x}}_\ast (X_0(x))),$$ that is, (see
(\ref{inv-izq})),
\begin{equation}\label{util'}
\rho _\ast (\phi _0)(x)=-\rho (X_0)(x).
\end{equation}
On the other hand, using (\ref{fi0}), (\ref{dif-induc}) and
(\ref{util}), it follows that $$e^{-\sigma} i(\delta \sigma )
(\izq{d_\ast X})=-i(\delta \sigma )({\cal L}_\ssli \Lambda
)+e^{-\sigma }(\phi _0(X)\circ \s)\izq{X_0}.$$ Consequently, using
again (\ref{fi0}), we have that
\begin{equation}\label{ecu5}
 i(\phi _0)(d_\ast X)=-i((\delta
\sigma )({\cal L}_\ssli \Lambda ))\circ \epsilon +\phi _0(X)X_0.
\end{equation}
Finally, from (\ref{compatibilidad}) and (\ref{fi0}), we deduce
that $$0= [\izq{X},\der{X_0}]=  i(\delta \sigma )({\cal L}_\ssli
\Lambda ) + \sostP (\delta (\phi _0(X)\circ \s)) - e^{-\sigma}
(\phi _0(X)\circ \s)\izq{X_0}+e^{-\sigma} \izq{\lcf X,X_0\rcf }$$
which implies that (see (\ref{diferencial}), (\ref{fi0}),
(\ref{3.8I}) and (\ref{ecu5})) $$i(\phi _0)(d_\ast X)+d_\ast(\phi
_0(X)) +\lcf X_0,X\rcf =0.$$\QED

Next, we will describe the generalized Lie bialgebroids associated
with some examples of Jacobi groupoids. We remark that two
generalized Lie bialgebroids $((A,\phi_0),(A^*,X_0))$ and
$((B,\omega_0),(B^\ast,Y_0))$ over a manifold $M$ are isomorphic
if there exists a  Lie algebroid isomorphism ${\cal I}:A\to B$
such that ${\cal I}(X_0)=Y_0$ and, in addition, the adjoint
operator ${\cal I}^\ast:B^\ast \to A^\ast $ is also a Lie
algebroid isomorphism satisfying ${\cal I}^\ast
(\omega_0)=\phi_0.$

\begin{examples}
{\rm {\bf 1.-}{\it Poisson groupoids}

 If $(G,\Jacobi ,\sigma )$
is a Jacobi groupoid with $E=0$ and $\sigma =0$, that is,
$(G,\Lambda )$ is a Poisson groupoid, then we have that $\phi _0$
and $X_0$ identically vanish (see (\ref{fi0}) and Remark
\ref{cocy0}). Therefore, (\ref{condcomp}) and Theorem \ref{bajada}
imply a well-known result (see \cite{MX}): if $(G,\Lambda )$ is a
Poisson groupoid then the pair $(AG,A^\ast G)$ is a Lie
bialgebroid.

{\bf 2.-}{\it Contact groupoids}

 Let $(\gr ,\eta ,\sigma )$ be a
contact groupoid and $(\Jacobi )$ be the Jacobi structure
associated with the contact 1-form $\eta$. Then, $(\gr ,\Jacobi
,\sigma )$ is a Jacobi groupoid.

Now, denote by $(\Lambda _0,E_0)$ the Jacobi structure on $M$
characterized by the conditions (\ref{Jacobi-base-contacto}), by
$X_0$ the section of the Lie algebroid $AG$ of $G$ satisfying
$E=-\der{X_0}$ and by ${\cal I}:T^\ast M\times\R\to AG$ the Lie
algebroid isomorphism given by (\ref{2.2'''}). If we consider the
section $(0,-1)\in \Omega ^1(M)\times C^\infty (M,\R)$ of the
vector bundle $T^\ast M\times\R\to M$, we have that (see
(\ref{2.2'''}))
\begin{equation}\label{3.7vii}
{\cal I}(0,-1)=X_0.
\end{equation}
Moreover, if ${\cal I}^\ast :A^\ast G\to TM\times\R$ is the
adjoint operator of ${\cal I}$, from (\ref{2.2'''}), it follows
that
\begin{equation}\label{3.7viii}
{\cal I}^\ast (\nu _x)=(-\s _\ast ^{\tilde{x}}(\sostP
(\tilde{\epsilon}(\nu _x))),-\nu _x(X_0(x)))
\end{equation}
for $\nu _x\in A^\ast _xG$, where $\tilde{\epsilon}$ is the
inclusion in the Lie groupoid $T^\ast G\gpd A^\ast G$.

Next, denote by $(\makebox{{\bf [}\, ,\, {\bf ]}}_-, \pi_-)$ the
Lie algebroid structure on the vector bundle $TM\times \R\to M$
defined by $$\begin{array}{l} \makebox{{\bf [}}
(X,f),(Y,g)\makebox{{\bf ]}}=(-[X,Y],-(X(g)-Y(f))),\qquad \pi
_-(X,f)=-X
\end{array}$$ for $(X,f)\in \mathfrak X(M)\times C^\infty (M,\R)$.

On the other hand, if on the vector bundle $TG\times\R\to G$ we
consider the natural Lie algebroid structure (see Section
\ref{algebroides}) then the map $\sostJ :T^\ast G\times\R\to
TG\times\R$ is a Lie algebroid homomorphism between the Lie
algebroids $(T^\ast G \times \R ,\lcf \, ,\, \rcf _{(\Lambda ,E)},
\widetilde{\#}_{(\Lambda ,E)})$ and $TG\times\R$. Using this fact,
(\ref{3.8I}) and since $M\cong \epsilon (M)$ is a coisotropic
submanifold of $G$, we deduce that ${\cal I}^\ast$ defines an
isomorphism between the Lie algebroids $A^\ast G$ and
$(TM\times\R, \makebox{{\bf [}\, ,\, {\bf ]}}_-, \pi_-)$. In
addition, from (\ref{3.7viii}) and Proposition \ref{utiles}, we
obtain that ${\cal I}^\ast (\phi _0)=(-E_0,0).$

In conclusion, if on the vector bundle $T^\ast M\times\R\to M$
(respectively, $TM\times\R\to M$) we consider the Lie algebroid
structure $(\lcf \, ,\, \rcf _{(\Lambda _0,E _0)},
\widetilde{\#}_{(\Lambda _0,E_0)})$ (respectively, $(\makebox{{\bf
[}\, ,\, {\bf ]}}_-, \pi_-)$) then the generalized Lie
bialgebroids $((AG,\phi_0),(A^\ast G,X_0))$ and $((T^\ast M\times
\R,(-E_0,0)),(TM\times \R,(0,-1))$ are isomorphic. Note that the
Jacobi structure on $M$ induced by the generalized Lie bialgebroid
$((T^*M\times \R,(-E_0,0)),(TM\times \R,(0,-1)))$ is just
$(\Lambda_0,E_0)$ (see (\ref{1.19I})).

{\bf 3.-}{\it Jacobi-Lie groups}

Let $G$ be a Lie group with identity element $\mathfrak e$,
$\mult$ be a multiplicative function and $(\Jacobi )$ be a Jacobi
structure on $G$ such that $\Lambda$ is $\sigma$-multiplicative,
$E$ is a right-invariant vector field and $$ \sostP(\delta \sigma
)_{(g)}=-E_g+\expsn (L_g)_\ast ^\mathfrak e (E(\mathfrak
e)),\mbox{ for all }g\in G.$$ Then, $(G\gpd \{\mathfrak e
\},\Jacobi ,\sigma )$ is a Jacobi groupoid (see Examples
\ref{ejemplos}).

The Lie algebroid of $G$ is just the Lie algebra $\mathfrak g$ of
$G$, that is, $AG=\mathfrak g$ and, from (\ref{fi0}), it follows
that $\phi _0=(\delta \sigma )(\mathfrak e)$.

On the other hand, since $\Lambda (\mathfrak e)=0$, one may
consider the intrinsic derivative $\delta _{\mathfrak
e}\Lambda:\mathfrak g\to \wedge ^2\mathfrak g$ of $\Lambda$ at
$\mathfrak e$. In fact, using (\ref{ecjacobi}) and (\ref{3.8I}),
we deduce that the Lie bracket $[\, ,\, ]_\ast$ on the dual space
$A^\ast G=\mathfrak g^\ast$ of $\mathfrak g$ is given by $$[\omega
,\nu ]_\ast =[\omega ,\nu ]_\Lambda -\omega (E(\mathfrak e))\nu
+\nu (E(\mathfrak e))\omega $$ for $\omega, \nu \in \mathfrak
g^\ast$, where $[\, ,\, ]_\Lambda:\mathfrak g^\ast \times\mathfrak
g^\ast\to \mathfrak g^\ast $ is the adjoint map of the intrinsic
derivative of $\Lambda$ at $\mathfrak e$. In addition, the
1-cocycle $X_0$ on $\mathfrak g^\ast$ is $X_0=-E(\mathfrak e)$.

Thus, using Theorem \ref{bajada}, we conclude that the pair
$((\mathfrak g, (\delta \sigma )(\mathfrak e)),(\mathfrak g^\ast
,-E(\mathfrak e)))$ is a generalized Lie bialgebroid over $\{
\mathfrak e\}$, i.e., a generalized Lie bialgebra. This result was
proved in \cite{IM2} (see Theorem 3.12 in \cite{IM2}).

{\bf 4.-}{\it An abelian Jacobi groupoid}

 Let $(L,\lcf \, ,\,\rcf
,\rho)$ be a Lie algebroid over a manifold $M$ and $\omega _0\in
\Gamma (L ^\ast)$ be a 1-cocycle of $L$. We may consider on $L^
\ast$ the Jacobi structure $(\Lambda_{(L^\ast ,\omega
_0)},E_{(L^\ast ,\omega _0)})$ given by (\ref{Jacobi-lineal}) and
the Lie groupoid structure for which $\s= \t$ is the vector bundle
projection $\tau :L^\ast \to M$ and the partial multiplication is
the addition in the fibers. As we know (see Examples
\ref{ejemplos}), $(L^\ast \gpd M,\Lambda_{(L^\ast ,\omega
_0)},E_{(L^\ast ,\omega _0)},0)$ is a Jacobi groupoid and we have
the corresponding generalized Lie bialgebroid $((A(L^\ast ),\phi
_0=0),(A^\ast(L ^\ast),X_0))$.

On the other hand, if $0:M\to L^\ast$ is the zero section of
$L^\ast$ and $\mu \in \tau \inv (x)=L_x^\ast$, we will denote by
$\mu ^v(0(x))\in T_{0(x)}L^\ast_x$ the vertical lift of $\mu$ to
$L^\ast$ at the point $0(x)$. Then, the map $$v:L^\ast \to
A(L^\ast ),\quad \mu \in L^\ast _x\mapsto \mu ^v(0(x))\in
A_x(L^\ast),$$ defines an isomorphism between the vector bundles
$L ^\ast$ and $A(L^\ast)$. Moreover, using (\ref{Jacobi-lineal})
and since $\s =\tau$ and the Lie bracket of two left-invariant
vector fields on $L^\ast$ is zero, we conclude that: i) $v$
defines an isomorphism between the Lie algebroid $L^\ast$ (with
the trivial Lie algebroid structure) and $A(L^\ast)$ and ii)
$v(\omega_0)=X_0$. In addition, if $v^\ast :A^\ast (L^\ast )\to L$
is the adjoint map of $v:L^\ast \to A(L^\ast )$ then, from
(\ref{ecjacobi}), (\ref{P-lineal}), (\ref{Jacobi-lineal}) and
(\ref{3.8I}), we deduce that $v^\ast$ induces an isomorphism
between the Lie algebroids $A^\ast (L^\ast )$ and $(L,\lcf \, ,\,
\rcf ,\rho )$.

Therefore, we have proved that the generalized Lie bialgebroids
$((A(L^\ast ),0),(A^\ast(L ^\ast),X_0))$ and  $((L^\ast
,0),(L,\omega _0))$ are isomorphic.

{\bf 5.-}{\it The banal Jacobi groupoid}

Let $(M,\Jacobi )$ be a Jacobi manifold and $G$ the product
manifold $M\times\R\times M$. Denote by $(\Lambda ',E ')$ the
Jacobi structure on $G$ given by (\ref{2.9''}) and by $\mult$ the
function defined by $\sigma (x,t,y)=t$. Then, one may consider a
Lie groupoid structure in $G$ over $M$ in such a way that $(\gr
,\Lambda ',E',\sigma )$ is a Jacobi groupoid (see Examples
\ref{ejemplos}). Thus, we have the corresponding generalized Lie
bialgebroid $((AG,\phi _0),(A^\ast G,X_0))$. As we know, the map
$\Phi :TM\times\R \to AG$ given by (\ref{2.9'}) defines an
isomorphism between the Lie algebroids $(TM\times \R
,\makebox{{\bf [}\, ,\, {\bf ]}}, \pi)$ and $AG$ and, moreover, it
follows that $\Phi (-E,0)=X_0$.

Now, let $\Phi ^\ast :A^\ast G\to T^\ast M\times\R$ be the adjoint
map of $\Phi$. Then, using (\ref{ecjacobi}), (\ref{2.9'}),
(\ref{2.9''}) and (\ref{fi0}), we deduce that $\Phi ^\ast$ induces
an isomorphism between the Lie algebroids $A^\ast G$ and $(T^\ast
M \times \R ,\lcf \, ,\, \rcf _{(\Lambda ,E)},
\widetilde{\#}_{(\Lambda ,E)})$ and, in addition, $\Phi ^\ast
(\phi _0)=(0,1)$.

Therefore, we have proved that the generalized Lie bialgebroids
$((AG,\phi _0),(A^\ast G,X_0))$ and  $((TM\times \R
,(0,1))$,$(T^\ast M \times \R ,(-E,0)))$ are isomorphic. }
\end{examples}
To finish this Section, we will relate the Jacobi structure on $G$
and the Jacobi structure on $M$ induced by the generalized Lie
bialgebroid structure of Theorem \ref{bajada}.
\begin{proposition}
Let $(\gr,\Jacobi ,\sigma )$ be a Jacobi groupoid and $(\Lambda
_0,E_0)$ be the Jacobi structure on $M$ induced by the generalized
Lie bialgebroid $((AG,\phi _0),(A^\ast G,X_0))$. Then, the
projection $\t$ is a Jacobi antimorphism between the Jacobi
manifolds $(G,\Jacobi)$ and $(M,\Lambda_0,E_0)$ and the pair $(\s
,e^\sigma)$ is a conformal Jacobi morphism.
\end{proposition}
\prueba Denote by $\{\;,\;\}$ (respectively, $\{\;,\;\}_0$) the
Jacobi bracket associated with the Jacobi structure $(\Jacobi)$
(respectively, $(\Lambda_0,E_0)$). Then, we must prove that
\[
\{\t^\ast f_1,\t^\ast f_2\}=-\t^\ast\{f_1,f_2\}_0,\;\;\;\;
e^{-\sigma}\{e^\sigma \s ^\ast f_1,e^\sigma \s ^\ast f_2\}=\s
^\ast\{f_1,f_2\}_0,
\]
for $f_1,f_2\in C^\infty(M,\R)$.

Now, if $(\rho_\ast,X_0):\Gamma(A^\ast G)\to {\mathfrak
X}(M)\times C^\infty(M,\R)$ is the map given by (\ref{3.8II}) and
$(\rho,\phi_0):\Gamma(AG)\to {\mathfrak X}(M)\times
C^\infty(M,\R)$ is the homomorphism of $C^\infty(M,\R)$-modules
defined by
\begin{equation}\label{3.16I}
(\rho,\phi_0)(X)=(\rho(X),\phi_0(X))
\end{equation}
then, from (\ref{1.19I}), (\ref{3.8II}) and (\ref{3.16I}), it
follows that
\begin{equation}\label{sostbase}
\# _{(\Lambda _0,E_0)}=(\rho _\ast ,X_0)\circ (\rho ,\phi _0)^\ast
,
\end{equation}
where $(\rho ,\phi _0)^\ast:\Omega^1(M)\times C^\infty(M,\R) \to
\Gamma(A^\ast G)$ is the adjoint operator of the homomorphism
$(\rho,\phi_0)$.

Using (\ref{TGR}) and since  $(\t ^T)_\sigma\circ \sostJ =(\rho
_\ast ,X_0)\circ \tilde{\t}_\sigma$ we have that
$$\begin{array}{lll}\{\t ^\ast f_1 ,\t ^\ast f_2\}&=&\langle
\sostJ (\t ^\ast \delta f_1,\t ^\ast f_1),(\t ^\ast \delta f_2,\t
^\ast f_2)\rangle \\ &=&\langle ((\t )^T_\sigma \circ \sostJ )(\t
^\ast \delta f_1,\t ^\ast f_1), (\delta f_2\circ\t,\t ^\ast f_2
)\rangle
\\ &=&\langle ((\rho _\ast ,X_0)\circ \tilde{\t}_\sigma )(\t
^\ast \delta f_1,\t ^\ast f_1), (\delta f_2\circ \t,\t ^\ast f_2
)\rangle .
\end{array}$$
From (\ref{inv-izq}), (\ref{T*G}), (\ref{T*GR}), (\ref{fi0}) and
(\ref{3.16I}), we deduce that $\tilde{\t}_\sigma ((\t
_\ast^g)^\ast( \omega_{\t(g)}) ,\lambda )=-(\rho ,\phi _0)^\ast$ $
(\omega_{\t(g)} ,\lambda )$, for $(\omega_{\t(g)} ,\lambda )\in
T^\ast_{\t(g)}M \times\R$. Using this fact and (\ref{sostbase}),
we get that $$\{\t ^\ast f_1 ,\t ^\ast f_2\}=\t ^\ast \{ f_1,f_2\}
_M.$$

On the other hand, using (\ref{inv-izq}), (\ref{1.19I}),
(\ref{TGR}), Proposition \ref{propiedades} and since $(\s
^T)_\sigma\circ \sostJ =(\rho _\ast ,X_0)\circ \tilde{\s}_\sigma$,
we obtain that $$\begin{array}{lll}e^{-\sigma}\{e^\sigma \s ^\ast
f_1 ,e^\sigma \s ^\ast f_2\}&=&e^{-\sigma} \langle \sostJ (\delta(
e^\sigma \s ^\ast f_1),e^\sigma\s ^\ast f_1),(\delta (e^\sigma \s
^\ast f_2),e^\sigma\s ^\ast f_2)\rangle
\\ &=&e^\sigma\langle ((\s ^T)_\sigma \circ \sostJ )(\s ^\ast \delta
f_1,\s ^\ast f_1), (\delta f_2\circ \s,\s ^\ast f_2 )\rangle
\\ && + e^\sigma (\s ^\ast f_1)
\langle\#_{(\Jacobi)}(\delta\sigma,0), (\s^\ast(\delta
f_2),\s^\ast f_2)\rangle\\ &=& e^\sigma \langle ((\rho _\ast
,X_0)\circ \tilde{\s}_\sigma )(\s ^\ast (\delta f_1),\s ^\ast
f_1), (\delta f_2\circ \s ,\s ^\ast f_2 )\rangle
\\ &&+ \s ^\ast( f_1E_0(f_2)).
\end{array}$$
Now, from (\ref{inv-izq}), (\ref{T*G}), (\ref{T*GR}), (\ref{fi0})
and  (\ref{3.16I}), it follows that
 $e^{\sigma(g)} \tilde{\s}_\sigma
((\s_\ast ^g) ^\ast (\omega_{\s (g)}) ,\lambda )=$ $(\rho ,\phi
_0)^\ast (\omega_{\s(g)} ,0 ),$ for $(\omega_{\s(g)} ,\lambda )\in
T^\ast_{\s(g)}M\times\R$. Therefore, $$e^{-\sigma}\{e^\sigma \s
^\ast f_1 ,e^\sigma \s ^\ast f_2\}=\s^\ast \{ f_1,f_2\}_0.$$

\vspace{-1cm}\QED

\subsection{Integration of generalized Lie bialgebroids}
In this Section, we will show a converse of Theorem \ref{bajada},
that is, we will show that one may integrate a generalized Lie
bialgebroid and obtain a Jacobi groupoid.
\subsubsection{Jacobi groupoids and Poisson groupoids} In this
first subsection, we will prove that a Poisson groupoid can be
obtained from any Jacobi groupoid and we will show the relation
between the generalized Lie bialgebroid associated with the Jacobi
groupoid and the Lie bialgebroid induced by the Poisson groupoid.

Let $\gr$ be a Lie groupoid and $\mult$ be a multiplicative
function. Then, using the multiplicative character of $\sigma$, we
can define a right action of $\gr$ on the canonical projection
$\pi _1:M\times \R \to M$ as follows
\begin{equation}\label{accion}
(x,t)\cdot g=(\s (g),\sigma (g)+t)
\end{equation}
for $(x,t)\in M\times \R$ and $g\in G$ such that $\t(g)=x.$ Thus,
we have the corresponding action groupoid $(M\times \R)\ast G\gpd
M\times \R.$ Moreover, if $(AG,\lcf\;,\;\rcf,\rho)$ is the Lie
algebroid of $G,$ the multiplicative function $\sigma$ induces a
$1$-cocycle $\phi_0$ on $AG$ given by (see (\ref{fi0}))
\begin{equation}\label{3.19}
\phi_0(x)(X_x)=X_x(\sigma),\;\;\; \mbox{for $x\in M$ and $X_x\in
A_x G.$}
\end{equation}
In addition, using the results in Section \ref{grupoides} (see
(\ref{1.9I})), we deduce that the $\R$-linear map
$\ast:\Gamma(AG)\to {\mathfrak X}(M\times \R)$ defined by
\begin{equation}\label{3.20}
X^\ast=(\rho(X)\circ \pi_1) + (\phi_0(X)\circ \pi_1)\frac{\partial
}{\partial t}
\end{equation}
induces an action of $AG$ on the projection $\pi_1:M\times \R\to
M$ and the Lie algebroid  of $(M\times \R)\ast G$ is just the
action Lie algebroid $AG\ltimes \pi_1$.

Now, it is easy to prove that $(M\times \R)\ast G$ may be
identified with the product manifold $G\times \R$ and, under this
identification, the structural functions of the Lie groupoid are
given by
\begin{equation}\label{3.21}
\begin{array}{rcll}
\s_\sigma(g,t)&=&(\s(g), \sigma(g) + t),&\mbox{ for } (g,t)\in
G\times \R,\\ \t_\sigma(h,s)&=&(\t(h),s),& \mbox{ for } (h,s)\in
G\times \R,\\ m_\sigma((g,t),(h,s)),&=&(gh,t),&\mbox{ if } \s
_\sigma(g,t)=\t _\sigma(h,s),\\ \epsilon_\sigma(x,t)&=&(\epsilon
(x), t),& \mbox{ for } (x,t)\in M\times \R,\\
\iota_\sigma(g,t)&=&(\iota(g),\sigma(g) + t),& \mbox{ for }
(g,t)\in G\times \R.
\end{array}
\end{equation}
Thus, if $A(G\times \R)$ is the Lie algebroid of $G\times \R$ and
$X\in A_{(x,t)}(G\times \R),$ it is clear that $X\in A_xG$ and
therefore the map
\begin{equation}\label{3.21I}
{\cal J}:A(G\times \R)\to AG\times \R,\;\; X\in A_{(x,t)}(G\times
\R)\to {\cal J}(X)=(X,t)\in A_xG\times \R
\end{equation}
defines an isomorphism of vector bundles. Furthermore, if on
$AG\times \R$ we consider the Lie algebroid structure
$(\lcf\;,\;\rcf^{\,\bar{ }\,\phi_0},\bar\rho^{\,\phi_0})$ given by
(\ref{corchbarra}) then ${\cal J}$ is a Lie algebroid isomorphism.
In conclusion, the Lie algebroid of the Lie groupoid $\grr$ may be
identified with $(AG\times
\R,\lcf\;,\;\rcf^{\,\bar{}\,\phi_0},\bar \rho^{\,\phi_0})$.

We also have the following result.

\begin{proposition}\label{p3.8}
Let $\gr$ be a Lie groupoid and $\sigma:G\to \R$ be a
multiplicative function. Suppose that $(\Jacobi)$ is a Jacobi
structure on $G$, that $\tilde{\Lambda}=e^{-t}(\Lambda +
\frac{\partial}{\partial t}\wedge E)$ is the Poissonization on
$G\times \R$ and that in $G\times \R$ we consider the Lie groupoid
structure on $M\times \R$ with structural functions given by
(\ref{3.21}). Then, $(\gr, \Jacobi,\sigma)$ is a Jacobi groupoid
if and only if $(\grr,\tilde{\Lambda})$ is a Poisson groupoid.
\end{proposition} \prueba From (\ref{TG}) and (\ref{3.21}),
it follows that the projections $(\s_\sigma)^T,(\t_\sigma)^T$, the
inclusion $(\epsilon_\sigma)^T$ and the partial multiplication
$\oplus_{T(G\times \Rp)}$ of the tangent groupoid $T(G\times
M)\gpd T(M\times \R)$ are given by
\begin{equation}\label{3.22}
\begin{array}{l}
\kern-16pt(\s_\sigma)^T(X_g + \lambda\frac{\partial}{\partial
t}_{|t}) \kern-3pt=\kern-3pt\s^T(X_g) + (\lambda +
X_g(\sigma))\frac{\partial }{\partial t}_{|t+\sigma(g)},\mbox{ for
} (X_g + \lambda\frac{\partial}{\partial t}_{|t})\kern-3pt\in
\kern-3ptT_{(g,t)} (G\times \R),\\ \kern-16pt(\t_\sigma)^T(Y_h +
\mu\frac{\partial}{\partial t}_{|s})= \t^T(Y_h) + \mu
\frac{\partial}{\partial t}_{|s},\mbox{ for } (Y_h +
\mu\frac{\partial}{\partial t}_{|s})\in T_{(h,s)} (G\times \R),
\\\kern-16pt (X_g +
\lambda\frac{\partial}{\partial t}_{|t})\oplus_{T(G\times \Rp)}(
Y_h + \mu\frac{\partial}{\partial t}_{|s})= X_g\oplus_{TG}Y_h +
\lambda\frac{\partial }{\partial t}_{|t}, \\ \kern-16pt
 (\epsilon_\sigma)^T(X_x
+ \lambda\frac{\partial}{\partial t}_{|t})=\epsilon^T(X_x) +
\lambda\frac{\partial}{\partial t}_{|t}, \mbox{ for } (X_x +
\lambda\displaystyle\frac{\partial }{\partial t}_{|t}) \in
T_{(x,t)}(M\times \R).
\end{array}
\end{equation}
On the other hand, using (\ref{T*G}) and (\ref{3.21}), we deduce
that the projections
$\widetilde{\s_\sigma},\widetilde{\t_\sigma}$, the inclusion
$\widetilde{\epsilon_\sigma}$ and the partial multiplication
$\oplus_{T^*(G\times \Rp)}$ in the cotangent groupoid $T^\ast
(G\times \R)\gpd A^\ast G\times \R$ are defined by
\begin{equation}
\begin{array}{l}\begin{array}{rcl}
\widetilde{\s_\sigma}(\omega_g + \gamma \delta
t_{|t})&=&(\tilde{\s}(\omega_g),\sigma(g) + t),\mbox{ for }
(\omega_g + \gamma \delta t_{|t})\in T^\ast_{(g,t)}(G\times \R)\\
\mbox{}\kern-30pt\widetilde{\t_\sigma}(\nu_h + \zeta \delta
t_{|s})&=&(\tilde{\t}(\nu_h)-\zeta
(\delta\sigma)_{\widetilde{\t(g)}} ,s),\mbox{ for } (\nu_h + \zeta
\delta t_{|s})\in T^\ast_{(h,s)}(G\times \R)\end{array}\\
\begin{array}{rcl} (\omega_g +\gamma \delta t_{|t})\oplus_{T^\ast(G\times
\Rp)}( \nu_h + \zeta \delta t_{|s})&=& (\omega_g +
\zeta(\delta\sigma)_g)\oplus_{T^\ast G}\nu_h + (\gamma +
\zeta)\delta t_{|t}\\
\widetilde{\epsilon_\sigma}(\omega_x,t)&=&\tilde\epsilon(\omega_x)
+ 0\,\delta t_{|t},\mbox{ for } (\omega_x,t)\in A_x^*G\times \R.
\end{array}\end{array}
\end{equation}
Moreover, from (\ref{homog}), we have that the homomorphism
$\#_{\tilde\Lambda}:T^\ast (G\times \R)\to T(G\times \R)$ is given
by
\begin{equation}\label{3.24}
\#_{\tilde\Lambda}(\omega_g + \gamma \delta
t_{|t})=e^{-t}(\#_\Lambda(\omega_g) + \gamma
E_g-\omega_g(E_g)\frac{\partial}{\partial t}_{|t}),
\end{equation}
for $(\omega_g + \gamma \delta t_{|t})\in T_{(g,t)}^\ast (G\times
\R).$

Now, we consider in $T^\ast G\times \R$ (respectively, $TG\times
\R$) the Lie groupoid structure over $A^\ast G$ (respectively,
$TM\times \R$) with structural functions defined by (\ref{T*GR})
(respectively, (\ref{TGR})). Then, an straightforward computation,
using (\ref{TG}), (\ref{T*G}), (\ref{TGR}), (\ref{morf-Jacobi}),
(\ref{T*GR}), (\ref{3.22})-(\ref{3.24}), shows that $\sostJ:T^\ast
G\times \R\to TG\times \R$ is a Lie groupoid morphism over some
map $\varphi_0:A^\ast G\to TM\times \R$ if and only if
$\#_{\tilde\Lambda}:T^\ast(G\times \R)\to T(G\times \R)$ is a Lie
groupoid morphism over some map $\tilde\varphi_0:A^\ast G\times
\R\to T(M\times \R).$ This proves the result. \QED

As we know (see Section \ref{glbi}) if $((A,\phi_0),(A^\ast,X_0))$
is a generalized Lie bialgebroid and on the vector bundle $A\times
M\to M\times \R$ (respectively, $A^\ast\times \R\to M\times \R)$
we consider the Lie algebroid structure $(\lcf\;,\;\rcf^{\,\bar{
}\,\phi_0},\bar\rho^{\,\phi_0})$ (respectively,
$(\lcf\;,\;\rcf^{\,\hat{ } \,X_0}_\ast,\hat{\rho}_\ast^{\, X_0}))$
 then the pair $(A\times \R,A^\ast \times \R)$ is a Lie bialgebroid. In
particular, if $(\gr,\Jacobi,\sigma)$ is a Jacobi groupoid and
$AG$ is the Lie algebroid of $G$ then the pair $(AG\times
\R,A^\ast G\times \R)$ is a Lie bialgebroid. Furthermore, we have

\begin{proposition}\label{Pois}
Let $(\gr ,\Jacobi ,\sigma )$ be a Jacobi groupoid and $(\grr
,\tilde{\Lambda })$ be the corresponding Poisson groupoid. If
$((AG,\phi _0),(A^\ast G,X_0))$ (respectively, $(A(G\times
\R),A^\ast (G\times \R))$) is the generalized Lie bialgebroid
(respectively, the Lie bialgebroid) associated with $(\gr ,\Jacobi
,\sigma )$ (respectively, $(\grr ,\tilde{\Lambda })$), then the
Lie bialgebroids $(A(G\times \R),A^\ast (G\times \R))$ and
$(AG\times \R, A^\ast G\times \R)$ are isomorphic.
\end{proposition}
\prueba Denote by $(\lcf\;,\;\rcf,\rho)$ the Lie algebroid
structure on $AG$ and by ${\cal J}:A(G\times \R)\to AG\times \R$
the isomorphism between the Lie algebroids $A(G\times \R)$ and
$(AG\times \R,\lcf\;,\;\rcf^{\,\bar{
}\,\phi_0},\bar{\rho}^{\,\phi_0})$ given by (\ref{3.21I}).

Now, let $\tilde{\cal J}^\ast:T^\ast G\times \R\times \R\to T^\ast
(G\times \R)$ be the map defined by
\[
\tilde{\cal J}^\ast (\omega_g,\gamma,t)=\omega_g+\gamma \delta
t_{|t},\mbox{ for }\omega_g\in T^*_gG \mbox{ and }\gamma ,t\in \R.
\]
Using the results in \cite{IM} (see Section 3.2 in \cite{IM}), we
deduce that
\begin{equation}\label{3.25}
\begin{array}{rcl}
\tilde{\cal J}^\ast \lcf
(\tilde\s,\tilde{f}),(\tilde{\t},\tilde{g})\rcf_{(\Jacobi)}^{\,\hat{
}\, X_0}&=&\lcf\tilde{\cal J}^\ast(\tilde\alpha,\tilde{f}),
\tilde{\cal
J}^\ast(\tilde\t,\tilde{g})\rcf_{\tilde\Lambda}=\lcf\tilde\s+
\tilde{f}\delta t, \tilde\t + \tilde{g}\delta
t\rcf_{\tilde\Lambda},\\ \#_{\tilde\Lambda}(\tilde{\cal J}^\ast
(\tilde\s,\tilde{f}))&=&\widehat{\tilde\sostJ}^{X_0}(\tilde\s,\tilde{f}),
\end{array}
\end{equation}
for $\tilde\s, \tilde\t$ time-dependent $1$-forms on $G$ and
$\tilde{f},\tilde{g}\in C^\infty(G\times \R,\R),$ where
$(\lcf\;,\;\rcf_{(\Jacobi)},\tilde\sostJ)$ (respectively,
$(\lcf\;,\;\rcf_{\tilde\Lambda},\#_{\tilde\Lambda})$) is the the
Lie algebroid structure on $T^\ast G\times \R$ (respectively,
$T^*(G\times \R))$ induced by the Jacobi structure $(\Jacobi)$
(respectively, the Poisson structure $\tilde\Lambda$) on $G$
(respectively, $G\times \R)$.

On the other hand, if we identify $A^\ast G$ (respectively,
$A^\ast (G\times \R)$) with the conormal bundle of $\epsilon (M)$
(respectively, $\epsilon _\sigma (M\times \R)$) then the
restriction of $\tilde{\cal J}^*$ to $A^\ast G\times \{0\}\times
\R\cong A^\ast G\times \R$ is just the adjoint operator ${\cal
J}^\ast:A^\ast G\times \R\to A^\ast (G\times \R)$ of ${\cal J}$.
Therefore, from (\ref{corchtilde}), (\ref{3.8I}), (\ref{3.25}) and
Remark \ref{cocy0}, we conclude that the map ${\cal J}^\ast$ is an
isomorphism between the Lie algebroids $(A^\ast G\times \R,
\lcf\;,\;\rcf_\ast^{\, \hat{ }\, X_0},\hat{\rho}_\ast^{\,X_0})$
and $ A^\ast(G\times \R).$ \QED
\subsubsection{Integration of generalized Lie bialgebroids} In this
Section, we will show a converse of Theorem \ref{bajada}.

For this purpose, we will use the notion of the derivative of an
affine $k$-vector field on a Lie groupoid (see \cite{MX2}). Let
$G$ be a Lie groupoid with Lie algebroid $AG$ and $P$ be an affine
$k$-vector field on $G.$ Then, the derivative of $P$, $\delta P$,
is the map $\delta P:\Gamma(AG)\to \Gamma(\wedge^k(AG))$ defined
as follows. If $X\in \Gamma(AG)$, $\delta P(X)$ is the element in
$\Gamma(\wedge^k(AG))$ whose left translation is ${\cal L}_{\ssli}
P.$

Now, we will prove the announced result at the beginning of this
Section.
\begin{theorem}\label{subida}
Let $((AG,\phi _0),(A^\ast G,X_0))$ be a generalized Lie
bialgebroid where $AG$ is the Lie algebroid of an $\s$-connected
and  $\s$-simply connected Lie groupoid $\gr$. Then, there is a
unique multiplicative function $\mult$ and a unique Jacobi
structure $(\Jacobi )$ on $G$ that makes $(\gr,\Jacobi,\sigma)$
into a Jacobi groupoid with generalized Lie bialgebroid $((AG,\phi
_0),(A^\ast G,X_0))$.
\end{theorem}
\prueba Since $G$ is $\alpha$-connected and $\alpha$-simply
connected, we deduce that there exists a unique multiplicative
function $\mult$ such that
\[
\phi_0(X)=X(\sigma),\makebox[1cm]{} \forall X\in \Gamma(AG).
\]
The multiplicative function $\mult$ allows us to construct a Lie
groupoid structure in $G\times \R$ over $M\times \R$ with
structural functions $\s_\sigma, \t_\sigma,
m_\sigma,\epsilon_\sigma$ and $\iota_\sigma$ given by
(\ref{3.21}).

If $(\lcf\;,\;\rcf,\rho)$ is the Lie algebroid structure on $AG$
then, as we know, the Lie algebroid of $G\times \R$ is $(AG\times
\R,\lcf\;,\;\rcf^{\,\bar{ }\,\phi_0},\bar\rho^{\,\phi_0}).$
Moreover, if $(\lcf\;,\;\rcf_\ast, \rho_\ast)$ is the Lie
algebroid structure on $A^\ast G$ and we consider on the vector
bundle $A^\ast G\times \R\to M\times \R$ the Lie algebroid
structure $(\lcf\;,\;\rcf_\ast^{\,\hat{ }\,
X_0},\hat\rho_\ast^{X_0})$ given by (\ref{corchtilde}), it follows
that the pair $(AG\times \R,A^\ast G\times \R)$ is a Lie
bialgebroid. Therefore, using Theorem 4.1 in \cite{MX2}, we obtain
that there is a unique Poisson structure $\tilde\Lambda$ on
$G\times \R$ that makes $G\times \R$ into a Poisson groupoid with
Lie bialgebroid $(AG\times \R, A^\ast G\times \R).$

We will see that the $2$-vector (on $G\times \R$) ${\cal
L}_{\frac{\partial }{\partial t}}\tilde\Lambda + \tilde{\Lambda}$
is affine. For this purpose, we will use the following relation
\begin{equation}\label{3.26}
{\cal L}_{\frac{\partial }{\partial t}} \izq{ {\tilde
P}}=\izq{\frac{\partial {\tilde{P}}}{\partial t}}, \mbox{ for }
\tilde{P}\in \Gamma(\wedge^k(AG\times \R)).
\end{equation}

Note that $\tilde{P}$ is a time-dependent section of the vector
bundle $\wedge^k(AG)\to M$ and, thus, one may consider the
derivative of $\tilde{P}$ with  respect to the time,
$\displaystyle\frac{\partial \tilde P}{\partial t}.$

From (\ref{3.26}) and Proposition \ref{sigma-afin}, we conclude
that the vector field $\displaystyle\frac{\partial }{\partial t}$
is affine. Consequently (see Proposition 2.5 in \cite{MX2}), the
$2$-vector ${\cal L}_{\frac{\partial }{\partial t}}\tilde\Lambda +
\tilde{\Lambda}$ is also affine.

Next, we will show that the Poisson structure $\tilde\Lambda$ is
homogeneous with respect to the vector field
$\displaystyle\frac{\partial }{\partial t}$. This fact implies
that $\tilde\Lambda$ is the Poissonization of a Jacobi structure
$(\Jacobi)$ on $G$ (see Remark \ref{poisonizacion}). Moreover,
from Propositions \ref{p3.8} and \ref{Pois}, we will have that
$(\gr, \Jacobi,\sigma)$ is a Jacobi groupoid with generalized Lie
bialgebroid $((AG,\phi_0),(A^\ast G,X_0)).$

Therefore, we must prove that $\tilde{\Lambda}$ is homogeneous.
Now, using Theorem 2.6 in \cite{MX2} and since $G$ is
$\alpha$-connected and the $2$-vector ${\cal
L}_{\frac{\partial}{\partial t}}\tilde\Lambda + \tilde\Lambda$ is
affine, we deduce that $\tilde{\Lambda}$ is homogeneous if and
only if:
\begin{enumerate}
\item The derivative of the $2$-vector ${\cal L}_{\frac{\partial
}{\partial t}}\tilde{\Lambda} + \tilde{\Lambda}$ is zero and
\item The restriction of the $2$-vector ${\cal L}_{\frac{\partial
}{\partial t}}\tilde{\Lambda} + \tilde{\Lambda}$ to the points of
$\epsilon_\sigma(M\times \R)$ is zero.
\end{enumerate}

First, we will show (i). If $H$ is a Poisson groupoid with Poisson
structure $\pi$ and Lie algebroid $AH$, we have that (see Theorem
3.1 in \cite{Xu})
\begin{equation}\label{3.27}
{\cal L}_\ssli\pi=-\izq{d_\ast X},
\end{equation}
for $X\in \Gamma(AH)$, where $d_*$ is the differential of the dual
Lie algebroid $A^\ast H$. Thus, from (\ref{3.26}) and
(\ref{3.27}), it follows that

$$
\begin{array}{ccl}
{\cal L}_\tssli \Big ({\cal L}_\fract \tLambda+\tilde{\Lambda}\Big
)&=& {\cal L}_{\fract }{\cal L}_\tssli \tLambda -{\cal
L}_{\tiny\izq{\frac{\partial \tilde{X}}{\partial t}}}\tLambda
+{\cal L}_\tssli \tLambda \\ &=&
\izq{\hat{d}_\ast^{X_0}\frac{\partial \tilde{X}}{\partial t}} -
\izq{\hat{d}_*^{X_0}\tilde{X}}-\izq{\frac{\partial (\hat{d}_\ast
^{X_0}\tilde{X})}{\partial t}},
\end{array}$$
for $\tilde{X}\in \Gamma(AG\times \R).$ On the other hand, using
the results in \cite{IM} (see Remark B.3 in \cite{IM}), we obtain
that $$\widehat{d_\ast}^{X_0} \tilde{Z}=e^{-t}\Big (
\hat{d}_{\ast}^{\, 0}\tilde{Z} +X_0\wedge (\tilde{Z} +
\frac{\partial \tilde{Z}}{\partial t})\Big ),\mbox{ for }
\tilde{Z}\in \Gamma(AG\times \R),$$ $\hat{d}_*^{\,0}$ being the
differential of the Lie algebroid $(A^\ast G\times
\R,\lcf\;,\;\rcf_\ast^{\,\hat{ }\, 0},\hat{\rho}_\ast ^{\, 0})$.
Consequently, we deduce that $${\cal L}_\tssli \Big ({\cal
L}_\fract \tLambda+\tilde{\Lambda}\Big )=0.$$ Next, we will show
(ii). If $(x,t)$ is a point of $M\times \R$ then $$ T^\ast
_{\epsilon_\sigma(x,t)}(G\times \R)\cong A^\ast_{(x,t)}(G\times
\R)\oplus ((\alpha_\sigma)_\ast^{(\tilde{x},t)})^\ast
(T^\ast_{(x,t)}(M\times \R)).$$ Therefore, it is enough to prove
that $$ ({\cal L}_\fract \tLambda+\tilde{\Lambda})(\delta
F_1,\delta F_2)_{|\epsilon_\sigma(M\times \Rp)}=0, $$ when $F_1$
and $F_2$ are either constant on $\epsilon_\sigma(M\times \R)$ or
equal to $(\alpha_\sigma)^\ast f_i,$ with $f_i\in C^\infty(M\times
\R,\R),$ $i=1,2$. We will distinguish three cases:

{\it First case}. Suppose that $F_1=(\alpha_\sigma)^*f_1$ and
$F_2=(\s_\sigma)^\ast  f_2,$ with $f_1,f_2\in C^\infty(M\times
\R,\R).$ Denote by $\tilde\Lambda_0$ the Poisson structure on
$M\times \R$ induced by the Lie bialgebroid $(AG\times \R,A^\ast
G\times \R)$ and by $\{\;,\;\}_{\tilde{\Lambda}}$ (respectively,
$\{\;,\;\}_{\tilde{\Lambda}_0}$) the Poisson bracket on $G\times
\R$ (respectively, $M\times \R$) associated with $\tilde{\Lambda}$
(respectively, $\tilde\Lambda_0$). Then, from Proposition
\ref{p3.8} and since the vector field
$\displaystyle\frac{\partial}{\partial t}$ on $G\times \R$ is
$\s_\sigma$-projectable, it follows that
\[\begin{array}{rcl}
({\cal L}_\fract \tLambda+\tilde{\Lambda})(\delta F_1,\delta
F_2)&=& \alpha_\sigma^\ast(\displaystyle\frac{\partial}{\partial
t}\{f_1,f_2\}_{\tilde\Lambda_0}-\displaystyle\{\frac{\partial
f_1}{\partial t},f_2\}_{\tilde\Lambda_0}
-\{f_1,\displaystyle\frac{\partial f_2}{\partial
t}\}_{\tilde\Lambda_0} + \{f_1,f_2\}_{\tilde\Lambda_0}).
\end{array}
\]
Thus, using that the Poisson structure $\tilde\Lambda_0$ is
homogeneous with respect to the vector field $\displaystyle
\frac{\partial }{\partial t}$ on $M\times \R$ (see Theorem
\ref{bialgebrizacion}), we obtain that
\[
 ({\cal L}_\fract
\tLambda+\tilde{\Lambda})(\delta F_1,\delta F_2)=0.
\]

\medskip

{\it Second case}. Suppose that $F_1=(\s_\sigma)^\ast f_1,$ with
$f_1\in C^\infty(M\times \R,\R)$ and that $F_2$ is constant on
$\epsilon_\sigma(M\times \R)$. Following the proof of Lemma 4.12
in \cite{MX2}, we deduce that
\begin{equation}\label{3.28}
\{(\alpha_\sigma)^\ast f,
H\}_{\tilde\Lambda}=\izq{((\hat{\rho}_\ast^{X_0})^\ast (\delta
f))}(H),
\end{equation}
for $f\in C^\infty(M\times \R,\R)$ and $H\in C^\infty(G\times
\R,\R).$ Note that (see (\ref{corchtilde}))
\begin{equation}\label{3.29}
(\hat{\rho}_\ast^{X_0})^\ast(\omega + g\delta
t)=e^{-t}((\rho_\ast)^\ast(\omega) + g X_0),
\end{equation}
for $g\in C^\infty(M\times \R,\R)$ and $\omega$ a time-dependent
$1$-form on $M.$ Therefore, from (\ref{3.26}), (\ref{3.28}),
(\ref{3.29}) and since $\displaystyle\frac{\partial F_2}{\partial
t}=0,$ we have that
\[
\begin{array}{rcl} ({\cal L}_\fract
\tLambda+\tilde{\Lambda})(\delta F_1,\delta
F_2)&\kern-10pt=&\kern-10pt -[\displaystyle\frac{\partial
}{\partial t},\izq{(\hat\rho_\ast ^{X_0})^\ast(\delta f_1)}](F_2)
- \izq{(\hat{\rho}_\ast^{X_0})^\ast(\delta(\frac{\partial
f_1}{\partial t}))}(F_2)- \izq{(\hat{\rho}_\ast^{X_0})^\ast(\delta
f_1)}(F_2)\\[5pt] &\kern-10pt=&\kern-10pt
\izq{\frac{\partial}{\partial t} ((\hat{\rho}_\ast^{X_0})^\ast
(\delta f_1))}-
\izq{(\hat{\rho}_\ast^{X_0})^\ast(\delta(\frac{\partial
f_1}{\partial t}))}- \izq{(\hat{\rho}_\ast^{X_0})^\ast(\delta
f_1)}(F_2) =0.
\end{array}
\]

\medskip

{\it Third case.} Suppose that $F_1$ and $F_2$ are constant on
$\epsilon_\sigma(M\times \R)$. Then, using that
$\epsilon_\sigma(M\times \R)$ is a coisotropic submanifold of
$(G\times \R,\tilde\Lambda),$ it follows that
\[
\{F_1,F_2\}_{\tilde\Lambda|\epsilon_\sigma(M\times \Rp)}=0.
\]
Moreover, since $\displaystyle\frac{\partial F_1}{\partial
t}=\displaystyle\frac{\partial F_2}{\partial t}=0$ and the
restriction to $\epsilon_\sigma(M\times \R)$ of the vector field
$\displaystyle\frac{\partial }{\partial t}$ is tangent to
$\epsilon_\sigma(M\times \R)$, we conclude that
\[
 ({\cal L}_\fract
\tLambda+\tilde{\Lambda})(\delta F_1, \delta F_2 )
_{|\epsilon_\sigma(M\times \Rp)}=0.
\]

\vspace{-25pt}

 \QED \vspace{-15pt}

\begin{examples}{\rm {\bf 1.-} {\em Lie bialgebroids}

Let $(AG,A^*G)$ be a Lie bialgebroid where $AG$ is the Lie
algebroid of an $\alpha$-connected and $\alpha$-simply connected
Lie groupoid $\gr$. Then, using Theorem \ref{subida} (see also
Examples \ref{ejemplos},1), we obtain that there exists a unique
Poisson structure $\Lambda$ on $G$ that makes $(\gr,\Lambda)$ into
a Poisson groupoid with Lie bialgebroid $(AG,A^*G)$. This result
was proved in \cite{MX2}.

{\bf 2.-} {\em Generalized Lie bialgebras}

If $G$ is a connected Lie group with identity element ${\mathfrak
e}$, $\mult$ is a multiplicative function and $\Lambda$ is a
$\sigma$-multiplicative $2$-vector such that the intrinsic
derivative of $\Lambda$ at ${\mathfrak e}$ is zero then $\Lambda$
identically vanishes (see \cite{IM2}).

Let $(({\mathfrak g},\phi_0)({\mathfrak g}^\ast, X_0))$ be a
generalized Lie bialgebra, that is, a generalized Lie bialgebroid
over a single point, and $G$ be a connected simply connected Lie
group with Lie algebra ${\mathfrak g}.$ Then, using (\ref{3.8I}),
Proposition \ref{propiedades} and Theorem \ref{subida} we deduce
the following facts: a) there exists a unique multiplicative
function $\mult$ and a unique $\sigma$-multiplicative $2$-vector
$\Lambda$ on $G$ such that $(\delta\sigma)({\mathfrak e})=\phi_0$
and the intrinsic derivative of $\Lambda$ at ${\mathfrak e}$ is
$-d_{\ast X_0},$ $d_{\ast X_0}$ being the $X_0$-differential of
the Lie algebra ${\mathfrak g}^\ast$; b) $\#_{\Lambda}(\delta
\sigma)=\der{X_0}- e^{-\sigma}\izq{X_0}$ and c) the pair
$(\Jacobi)$ is a Jacobi structure on $G$, where $E=-\der{X_0}.$
These results were proved in \cite{IM2} (see Theorem 3.10 in
\cite{IM2}). }
\end{examples}

\vspace{-.2cm} {\small {\bf Acknowledgments.} Research partially
supported by DGICYT grant BFM 2000-0808. D. Iglesias-Ponte wishes
to thank the Spanish Ministry of Education and Culture for an FPU
grant.}

\end{document}